\documentclass[a4paper,11pt]{article}
\usepackage[utf8]{inputenc}
\usepackage{url}
\usepackage{a4wide}
\usepackage{rotating}
\usepackage{lscape}
 \usepackage{xcolor}
\usepackage{algorithmic}
\usepackage{multirow}
\usepackage{array}
\usepackage{booktabs}
\usepackage{ragged2e}
\usepackage{makecell}
\usepackage{adjustbox}
\usepackage[normalem]{ulem}

\usepackage[driverfallback=dvipdfm]{hyperref}

\usepackage{setspace}
\usepackage{fancyhdr}
\usepackage{graphics}
\usepackage{textcomp}
\usepackage{amssymb}
\usepackage{amsmath}
\usepackage{cancel}

\usepackage[subrefformat=parens,labelformat=parens]{subcaption}
\usepackage[labelformat=simple,labelsep=quad,skip=3pt]{caption}
\usepackage{graphicx}

\newcounter{AlgCount}
\addtocounter{AlgCount}{-1}
\newcommand{\ralg}[1]{\refstepcounter{AlgCount}\label{#1}}

\usepackage{color}

\newcommand{\algrule}[1][.2pt]{\par\vskip.5\baselineskip\hrule height #1\par\vskip.5\baselineskip}

\def\R{{\mathbb{R}}}

\def\Tr{{\mathrm{Trace}}}

\def\rr{\right}
\def\l{\left}
\def\bc{\begin{center}}
\def\ec{\end{center}}
\def\bt{\begin{tabular}}
\def\et{\end{tabular}}

\def\norm#1{\l\| #1 \rr\| }
\def\br#1{\l(#1\rr)}

\def\ba{\begin{array}}
\def\ea{\end{array}}
\def\beq{\begin{equation}}
\def\eeq{\end{equation}}
\def\beqnn{\begin{equation*}}
\def\eeqnn{\end{equation*}}
\def\beann{\begin{eqnarray*} }
\def\eeann{\end{eqnarray*}}
\def\bea{\begin{eqnarray}}
\def\eea{\end{eqnarray}}
\def\bmi{\begin{minipage}[t]}
\def\emi{\end{minipage}}

\def\norm#1{\l\| #1 \rr\| }

\title{\bf{A Block Coordinate Descent-based Projected Gradient Algorithm for Orthogonal Non-negative Matrix Factorization
 }}

\author{Soodabeh Asadi\thanks{Institute for Data Science, School of Engineering, University of Applied Sciences and Arts, Northwestern Switzerland. Email:~soodabeh.asadidezaki@fhnw.ch}
 \and Janez Povh\thanks{University of Ljubljana, Faculty of Mechanical Engineering, A\" sker\" ceva ulica 6, 1000 Ljubljana, Slovenia. Email: janez.povh@fs.uni-lj.si}}

\begin{document}

\maketitle

\begin{abstract}
This article utilizes the projected gradient method (PG) for a non-negative matrix factorization problem (NMF), where one or both matrix factors must have orthonormal columns or rows. We penalise the  orthonormality constraints and apply the PG method via a block coordinate descent approach. This means that at a certain time one matrix factor is fixed and the other is updated by moving along the steepest descent direction computed from the penalised objective function and projecting onto the space of non-negative matrices.

Our method is tested on two sets of synthetic data for various values of penalty parameters. The performance is compared to the well-known multiplicative update (MU) method from Ding (2006), and with a modified global convergent variant of the MU algorithm recently proposed by Mirzal (2014). We provide extensive numerical results coupled with appropriate visualizations, which demonstrate that our method is very competitive and usually   outperforms the other two methods.

\end{abstract}
{\bf{Keywords}}: Non-negative matrix factorization, Orthogonality conditions, Projected gradient method, Multiplicative update algorithm, Block coordinate descent  \\

\section{Introduction}
\subsection{Motivation}
Many machine learning applications require processing large and high dimensional data. The data could be images, videos, kernel matrices, spectral graphs, etc., represented as an $m\times n$ matrix $R$. The data size and the amount of redundancy increase rapidly when $m$ and $n$ grow. To make the analysis and the interpretation easier, it is favorable to obtain compact and concise low rank approximation of the original data $R$.
This low-rank approximation is known to be very efficient in a wide range of applications, such as: text mining \cite{Berry-et-al2007,Pauca-et-al,Shahnaz-et-al}, document classification \cite{Berry-et-al2}, clustering \cite{Li-Ding,Wei-Liu-Gong}, spectral data analysis \cite{Berry-et-al2007,Kaarna}, face recognition \cite{Zafeiriou-et-al}, and many more.\\[2mm]

There exist many different low rank approximation methods. For instance, 
two well-known strategies, broadly used for data analysis, are singular value decomposition (SVD) \cite{Golub-Reinsch-1970} and
principle component analysis (PCA) \cite{Jolliffe-2005}. Much of real-world data are non-negative, and the related hidden parts express physical features only when the non-negativity holds. The factorizing matrices in SVD or PCA can have negative entries, making it hard or impossible to put a physical interpretation on them.
Non-negative matrix factorization was introduced as an attempt to overcome this drawback, i.e., to provide the desired low rank non-negative matrix factors. 
\subsection{Problem formulation}
A non-negative matrix factorization problem (NMF) is a problem of factorizing the input non-negative  matrix
 $R$ into the product of two lower rank non-negative matrices $G$ and $H$:
\bea \label{def NMF}
R\approx GH,
\eea
where $R\in \R_+^{m\times n}$ usually corresponds to the data matrix, $G\in \R_+^{m\times p}$ represents the basis matrix, and $H\in \R_+^{p\times n}$ is the coefficient matrix. With $p$ we denote the number of factors for which it is desired that $p\ll \min (m,n)$. If we consider each of the $n$ columns of $R$ being a sample of $m$-dimensional vector data, the factorization  represents each instance (column) as a non-negative linear combination 
of the columns of $G$, where the coefficients correspond to the columns of $H$. The columns of $G$ can be therefore interpreted as the $p$ pieces that constitute the data $R$.
To compute $G$ and $H$, condition \eqref{def NMF} is usually rewritten as a minimization problem using the Frobenius norm:
\begin{equation} \label{NMF}\tag{NMF}
\min_{G,H} f(G,H)=\frac12\norm{R-GH}_F^2,~~~G\geq 0,~H\geq 0.
\end{equation}

It is demonstrated in certain applications that 
the performance of the standard NMF in \eqref{NMF} can often be improved by adding auxiliary constraints which could be sparseness, smoothness, and orthogonality. Orthogonal NMF (ONMF) was introduced by Ding et al., \cite{Ding-Li-Peng-Park}. To improve the clustering capability of the standard NMF, they imposed orthogonality constraints on columns of $G$ or on rows of $H$. Considering the orthogonality on columns of $G$, it is formulated as follows:
\begin{equation}\label{ONMF} \tag{ONMF}
\min_{G,H} f(G,H)=\frac12\norm{R-GH}_F^2,~\mbox{s.t.}~G\geq 0,~H\geq 0, ~G^TG=I.
\end{equation}
If we  enforce orthogonality on the columns of $G$ and on  rows of $H$, we obtain the bi-orthogonal ONMF \eqref{biONMF}, which 
is formulated as 
\begin{equation}\label{biONMF} \tag{bi-ONMF}
\min_{G,H} f(G,H)=\frac12\norm{R-GH}_F^2,~\mbox{s.t.}~G\geq 0,~H\geq 0, ~G^TG=I, ~HH^T=I,
\end{equation}

where $I$ denotes the identity matrix. 

\subsection{Related work}

The NMF was firstly studied by Paatero et al., \cite{Paatero-Tapper,Anttila-Paatero-UntoTapper-Jarvinen} and was made popular by Lee and Seung \cite{Lee-Seung,Lee-Seung2001}.
There are several different existing methods to solve \eqref{NMF}. The most used approach to minimize \eqref{NMF} is a simple MU method proposed by Lee and Seung \cite{Lee-Seung,Lee-Seung2001}. In Chu et al., \cite{Chu-Diele}, several gradient-type approaches have been mentioned.
Chu et al., reformulated \eqref{NMF} as an unconstrained optimization problem, and then applied the standard gradient descent method. 
Considering both $G$ and $H$ as variables in \eqref{NMF}, it is obvious that $f(G,H)$ is a non-convex function. However, considering $G$ and $H$ separately, we can find two convex sub-problems. Accordingly, a block-coordinate descent (BCD) approach \cite{Lee-Seung2001} is applied to obtain values for $G$ and $H$ that correspond to a local minimum of $f(G,H)$.
Generally, the scheme adopted by BCD algorithms is to recurrently update blocks of variables only, while the remaining variables are fixed.
NMF methods which adopt this optimization technique are, e.g., the MU rule \cite{Lee-Seung}, the active-set-like method \cite{Kim-Park}, or the PG method for NMF \cite{Lin2007}. 
In \cite{Lin2007}, two PG methods were proposed for the standard NMF. The first one is an alternating least squares (ALS) method using projected gradients. This way, $H$ is fixed first and a new $G$ is obtained by PG. Then, with the fixed $G$ at the new value, the PG method looks for a new $H$. The objective function in each least squares problem is quadratic. This enabled the author to use Taylor's extension of the objective function to obtain an equivalent condition with the Armijo rule, while checking the sufficient decrease of the objective function as a termination criterion in a step-size selection procedure. The other method proposed in \cite{Lin2007} is a direct application of the PG method to \eqref{NMF}.
There is also a hierarchical ALS method for NMF which was originally proposed in  \cite{Cichocki-Zdunek-Amari,Halko-Martinsson-Tropp}
as an improvement to the ALS method. It consists of a BCD method with single component vectors as coordinate blocks.

As the original ONMF algorithms in \cite{Li-Ding,Wei-Liu-Gong} and their variants \cite{Yoo-Choi2008,Yoo-Choi2010,Choi} are all based on the MU rule, there has been no convergence guarantee for these algorithms. For example, Ding et al., \cite{Ding-Li-Peng-Park} only prove that the successive updates of the orthogonal factors will converge to a local minimum of the problem. Because the orthogonality constraints cannot be rewritten into a non-negatively constrained ALS framework, convergent algorithms for the standard NMF (e.g., see \cite{Lin2007,Kim-Sra-Dhillon2008,Kim-Sra-Dhillon2007,Kim-Park2008}) cannot be used for solving the ONMF problems. 
Thus, 
no convergent algorithm was available for ONMF until recently. Mirzal \cite{Mirzal2014} developed a convergent algorithm for ONMF. The proposed algorithm was designed by generalizing the work of Lin \cite{Lin2007-IEEETransactions} in which a convergent algorithm was provided for the standard NMF based on a modified version of the additive update (AU) technique of Lee \cite{Lee-Seung2001}. Mirzal \cite{Mirzal2014} provides the global convergence for his algorithm solving the ONMF problem. In fact, he first proves the non-increasing property of the objective function evaluated by the sequence of the iterates. Secondly, he shows that every limit point of the generated sequence is a stationary point, and finally he proves that the sequence of the iterates possesses a limit point. 
 \subsection{Our contribution}
In this paper, we consider the penalty reformulation of \eqref{biONMF}, i.e., we add the orthogonality constraints multiplied with penalty parameters to the objective function to obtain reformulated problems \eqref{ONMF} and \eqref{biONMF}. The main contributions are:
\begin{itemize}
    \item We develop an algorithm for \eqref{ONMF} and \eqref{biONMF}, which is essentially a BCD algorithm, in literature also known as alternating minimization, coordinate relaxation, the Gauss-Seidel method,
subspace correction, domain decomposition, etc., see e.g. \cite{Bertsekas-1999,Richarik_Takacs:14}. For each block optimization, we use a PG method and   Armijo rule to find a suitable step-size. 
\item We construct  synthetic data sets of instances for \eqref{ONMF} and \eqref{biONMF}, for which we know the optimum value by construction.
\item We use  MATLAB \cite{matlab} to implement our algorithm  and  two well-known  (MU-based) algorithms: the algorithm of Ding \cite{Ding-Li-Peng-Park} and of Mirzal 
\cite{Mirzal2014}. The code is available upon request.
\item The implemented algorithms are compared on the constructed synthetic data-sets in terms of: (i) the accuracy of the reconstruction, and (ii) the deviation of the factors from orthonormality. 
Accuracy is measured by the  so-called root-square error (RSE), defined as 
\bea \label{relative error}
\text{RSE}:=\frac{\norm{R-GH}_F}{1+\norm{R}_F},
\eea
and deviations from orthonormality are computed using formulas \eqref{def:vio} and \eqref{def:vio2} from  Sect.~\ref{sec:num}.
Our numerical results show that our algorithm is very competitive and almost always outperforms the MU algorithms.
\end{itemize}
\subsection{Notations}
Some notations used throughout our work are described here. We denote scalars and indices by lower-case Latin letters, vectors by lowercase boldface Latin letters, and matrices by capital Latin letters. $\R^{m\times n}$ denotes the set of $m$ by $n$ real matrices, and $I$ symbolizes the identity matrix. We use the notation $\nabla$ to show the gradient of a real-valued function. We define $\nabla^+$ and $\nabla^-$ as the positive and (unsigned) negative parts of $\nabla$, respectively, i.e., $\triangledown=\triangledown^+-\triangledown^-$.  $\odot$ and $\oslash$ denote the element-wise multiplication and the element-wise division, respectively.
\subsection{Structure of the paper}
 The rest of our work is organized as follows. In Sect. 2, we review the well-known MU method and the rules being used for updating the factors per iteration in our computations. We also outline the global convergent MU version of Mirzal \cite{Mirzal2014}. We then present our PG method and discuss the stopping criteria for it. Sect. 4 presents the synthetic data and the result of implementation of the three decomposition methods presented in Sect. 3. This implementation is done for both the problem \eqref{ONMF}, as well as \eqref{biONMF}.  Some concluding results are presented in Sect. 5.

\section{Existing methods to solve \eqref{NMF}}
\subsection{MU method of Ding \cite{Ding-Li-Peng-Park}}\label{subsec:DING}
Several popular approaches to solve 
\eqref{NMF} are based on so-called MU algorithms, which are simple to implement and often yield good results. The MU algorithms originate from the work of Lee and Seung \cite{Lee-Seung2001}. Various MU variants were later proposed by several researchers,  for an overview see \cite{Cichocki:09}. At each iteration of these methods, the elements of $G$ and $H$ are multiplied by certain updating factors.

As already mentioned, \eqref{ONMF} was proposed by Ding et al., \cite{Ding-Li-Peng-Park} as a tool to improve the clustering capability of the associated optimization approaches.
To adapt the MU algorithm for this problem, they employed standard Lagrangian techniques: they
introduced the Lagrangian multiplier $\Lambda$ (a symmetric matrix of size $p \times p$) for the orthogonality constraint, and minimized the Lagrangian function where the orthogonality constraint is moved to the objective function as the penalty term $\Tr(\Lambda (G^TG-I))$. 
The complementarity conditions from the related KKT conditions can be rewritten as a fixed point relation, which finally can lead to the following MU rule for \eqref{ONMF}:
\beq \label{Ding G-update-1}
\begin{array}{ccc}
G_{ij}=G_{ij}\sqrt{\frac{(RH^T)_{ij}}{(GG^TRH^T)_{ij}}},~i=1, \cdots, m,~j=1,\cdots ,p,\\
H_{st}=H_{st}\sqrt{\frac{(R^TG)_{st}}{(H^TG^TG)_{st}}},~s=1,\cdots, p,~t=1,\cdots, n.
\end{array}
\eeq
They extended this approach to 
 non-negative three factor factorization with demand that two factors satisfy orthogonality conditions, 
 which is a generalization of
\eqref{biONMF}. 
The MU rules (28)-(30) from \cite{Ding-Li-Peng-Park}, adapted  to  \eqref{biONMF}, are the main ingredients of Algorithm \ref{alg:Ding}, which we will call Ding's algorithm. 
\begin{figure} [!htb]\ralg{alg:Ding0}
\algrule[1pt]
  {\bf{Algorithm 1. Ding's MU algorithm for \eqref{biONMF}}} \vspace{-3mm}\\
\algrule[1pt]
{\bf INPUT:} $R\in \R_+^{m\times n}$, $p\in N$
\begin{enumerate}
\item {\bf Initialize:} generate $G\geq 0$ as an $m \times p$ random matrix and  $H\geq 0$ as a $p\times n$ random matrix.\\
\item {\bf Repeat}
\beq \label{Ding G-update}
\begin{array}{ccc}
G_{ij}=G_{ij}\sqrt{\frac{(RH^T)_{ij}}{(GG^TRH^T)_{ij}+\delta}},~i=1, \cdots, m,~j=1,\cdots ,p,\\
H_{st}=H_{st}\sqrt{\frac{(G^TR)_{st}}{(G^TRH^TH)_{st}+\delta}},~s=1,\cdots, p,~t=1,\cdots, n.
\end{array}
\eeq
\item \textbf{Until} convergence or a maximum number of iterations or maximum time is reached.
\end{enumerate}
\textbf{OUTPUT:} $G,H$.
\algrule[1pt]
\vspace{-0mm}\ralg{alg:Ding}
\end{figure}
%
Algorithm \ref{alg:Ding} converges in the sense that the solution pairs $G$ and $H$ generated by this algorithm yield a sequence of decreasing $\text{RSE}$s, see 
\cite[Theorems 5,~7]{Ding-Li-Peng-Park}. 

If $R$ has zero vector as columns or rows, a division by zero may occur. In contrast, denominators close to zero may still cause numerical problems. To escape this situation, we follow  \cite{Piper-Pauca-Plemmons-Maile} and  add a small positive number $\delta$
to the denominators of the MU terms \eqref{Ding G-update}.
Note that Algorithm \ref{alg:Ding} can be easily adapted to solve \eqref{ONMF} by replacing the second MU rule from \eqref{Ding G-update} with the second MU rule of \eqref{Ding G-update-1}.
\subsection{MU method of Mirzal \cite{Mirzal2014} } \label{subsec:Mirzal}
 In \cite{Mirzal2014}, Mirzal proposed an algorithm for
 \eqref{ONMF} which is designed by generalizing the work of Lin \cite{Lin2007-IEEETransactions}.
 Mirzal used the so-called modified 
 additive update rule (the MAU rule), where the updated term is added to the current value for each of the factors. This additive rule has been used by Lin in  \cite{Lin2007-IEEETransactions} in the context of a
 standard NMF.   
 He also provided in his paper a convergence proof, stating that 
 the iterates generated by his algorithm converge in the sense that $\text{RSE}$ is decreasing and the limit point is a stationary point.
   In \cite{Mirzal2014}, Mirzal discussed the
orthogonality constraint on the rows of $H$, while in \cite{MirzalUnpublished} the same results are developed for the case of \eqref{biONMF}.

Here we review the Mirzal's algorithm for \eqref{biONMF}, presented in the unpublished paper \cite{MirzalUnpublished}. This algorithm actually solves the equivalent problem \eqref{ONMF-M} 
where the orthogonality constraints are moved into the objective function (the so-called penalty approach), and the importance of the orthogonality constraints are controlled by the penalty parameters $\alpha,\beta$: 
\begin{equation}
\label{ONMF-M} \tag{pen-ONMF}
\begin{array}{rcl}
\min_{G,H} F(G,H)&=&\frac12\norm{R-GH}_F^2+\frac{\alpha}{2}\norm{HH^T-I}_F^2+\frac{\beta}{2}\norm{G^TG-I}_F^2,\\[2mm]
\text{s.t.}& & G\geq0, ~H\geq 0
\end{array}
\end{equation} 
The gradients of the objective function with respect to $G$ and $H$ are:
\beq \label{Mirzal gradeint with respect to G}
\begin{array}{ccc}
\nabla_G f(G,H)=GHH^T-RH^T+\beta GG^TG-\beta G,\\
\nabla_H f(G,H)=G^TGH-G^TR+\alpha HH^TH-\alpha H.
\end{array}\eeq
For the objective function in \eqref{ONMF-M}, Mirzal proposed the MAU rules along with the use of $\bar{G}=(\bar{g})_{ij}$ and $\bar{H}=(\bar{h})_{ij}$, instead of $G$ and $H$, to avoid the zero locking phenomenon \cite[Section 2]{Mirzal2014}:
\bea\label{def of Gbar}
\bar{g}_{ij}=\begin{cases}
g_{ij},     &\mbox{if}~\nabla_Gf(G,H)_{ij} \geq0\\
\max\{g_{ij},\nu\},    &\mbox{if}~\nabla_Gf(G,H)_{ij} <0\\
\end{cases}\\
\label{def of Hbar}
\bar{h}_{st}=\begin{cases}
h_{st},     &\mbox{if}~\nabla_Hf(G,H)_{st} \geq0\\
\max\{h_{st},\nu\},    &\mbox{if}~\nabla_Hf(G,H)_{st} <0\\
\end{cases}
\eea
where $\nu$ is a small positive number.

Note that, the algorithms working with the MU rules for \eqref{ONMF-M} must be initialized with positive matrices to avoid zero locking from the start, but non-negative matrices can be used to initialize the algorithm working with the MAU rules (see \cite{MirzalUnpublished}).

Mirzal \cite{MirzalUnpublished} used the MAU rules with some modifications by considering $\bar{G}$ and $\bar{H}$ in order to guarantee the non-increasing property,
with a constant step to make $\delta_G$ and $\delta_H$ grow in order to satisfy the property.
Here, $\delta_G$ and $\delta_H$ are the values added within the MAU terms 
to the denominator of update terms for $G$ and $H$, respectively.
The proposed algorithm by Mirzal \cite{MirzalUnpublished} is summarised as Algorithm \ref{alg:Mirzal} below.
\begin{figure} [!htb]
\algrule[1pt]
   {\bf{Algorithm 2. Mirzal's algorithm for bi-ONMF} \cite{MirzalUnpublished}} \vspace{-3mm}\\
\algrule[1pt]
{\bf INPUT:} inner dimension $p$, maximum number of iterations: maxit; small positive $\delta$,  small positive $step$ to increase $\delta$.
\begin{enumerate}
    \item 
Compute initial $G^0\geq 0$ and  $H^0\geq 0$. 
\item {\bf For} $k=0: \text{maxit}$\\
$~~~~~~~\delta_G=\delta$;\\
$~~~~~~~\text{{\bf Repeat}}$\\
$~~~~~~~~~~~~g_{ij}^{(k+1)}=g_{ij}^{(k)}-\frac{\bar{g}_{ij}^{(k)}\times\nabla_Gf(G^{(k)},H^{(k)})_{ij}}{(\bar{G}^{(k)}H^{(k)}H^{(kT)}+\beta\bar{G}^{(k)}\bar{G}^{(kT)}\bar{G}^{(k)})_{ij}+\delta_G^{(k)}},~~~i=1\cdots m, ~j=1,\cdots p$;\\
$~~~~~~~~~~~\delta_G=\delta_G\times \text{step}$;\\
$~~~~~~~\text{{\bf Until}}~f(G^{(k+1)},H^{(k)})\leq f(G^{(k)},H^{(k)})$\\
$~~~~~~~\delta_H=\delta$;\\
$~~~~~~~\text{{\bf Repeat}}$\\
$~~~~~~~h_{st}^{(k+1)}=h_{st}^{(k)}-\frac{\bar{h}^{(k)}_{st}\times\nabla_Hf(G^{(k+1)},H^{(k)}){st}}{(G^{(k+1)T}G^{(k+1)}\bar{H}^{(k)}+\alpha \bar{H}^{(k)}\bar{H}^{(kT)}\bar{H}^{(k)})_{st}+\delta_H^{(k)}}, ~~~s=1,\cdots p,~t=1,\cdots n$;
$~~~~~~~~~~~\delta_H=\delta_H\times \text{step}$;\\
$~~~~~~~\text{{\bf Until}}~f(G^{(k+1)},H^{(k+1)})\leq f(G^{(k+1)},H^{(k)})$\\
$~~~~~~~\delta_H=\delta$;
\end{enumerate}
{\bf OUTPUT:} $G, H$.
\vspace{-1mm}
\algrule[1pt]
\vspace{-0mm}\ralg{alg:Mirzal}
\end{figure}

\section{PG method for \eqref{ONMF} and \eqref{biONMF}}
\label{sec:PGM}
\subsection{Main steps of PG method}
In this subsection we adapt the PG method   proposed by Lin
\cite{Lin2007} to solve both \eqref{ONMF} as well as  \eqref{biONMF}.
Lin applied PG to \eqref{NMF} in two ways. The first approach
is actually a BCD method.
This method consecutively fixes one block of variables ($G$ or $H$) and minimizes the simplified problem in the other variable.
The second approach by Lin directly minimizes \eqref{NMF}.
Lin's main focus was on the first approach and we follow it.
We again try to solve the penalised version of the problem \eqref{ONMF-M} by the block coordinate  descent method, which is summarised in Algorithm \ref{alg:PG1}.

\begin{figure} [!htb]
\algrule[1pt]\ralg{alg:PG1}
  {\bf{Algorithm 3. BCD method for \eqref{ONMF-M}}} \vspace{-3mm}\\
 \algrule[1pt]
{\bf  INPUT:} inner dimension $p$, initial matrices $G^0,~H^0$.
\begin{enumerate}
\item [1.] Set $k=0$.
\item [2.] \textbf{Repeat}
\begin{itemize}
\item[] 
 Fix $H:=H^k$ and  compute new $G$ as follows:
\bea \label{subprobG}
G^{k+1}:=\text{argmin}_{G\geq0} 
\frac12\norm{R-GH^k}_F^2 + \frac{\alpha}{2}\norm{H^kH^{kT}-I}_F^2+\frac{\beta}{2}\norm{G^TG-I}_F^2
\eea
\item[]  Fix $G:=G^{k+1}$ and  compute new $H$ as follows:
\bea \label{subprobH}
H^{k+1}\!:=\text{argmin}_{H\ge 0} 
\frac12\!\norm{R-G^{k+1}H}_F^2+\!\!\frac{\alpha}{2}\!\norm{HH^T\!-I}_F^T +\frac{\beta}{2}\!\norm{G^{(k+1)T}G^{k+1}\!-I}_F^2
\eea
\item[] $k:=k+1$
\end{itemize}
\item[3.] \textbf{Until} some stopping criteria is satisfied
\end{enumerate}
{\bf OUTPUT:} $G,H$.
\algrule[1pt]
\end{figure}

The objective function in \eqref{ONMF-M} is not quadratic any more, so we lose the nice properties about Armijo's rule that represent advantages for Lin. We managed to use the Armijo rule directly and still obtained good numerical results, see Sect. \ref{sec:num}.

 We refer to \eqref{subprobG} or \eqref{subprobH} as  sub-problems. Obviously, solving these sub-problems  in every iteration could be more costly than  Algorithms \ref{alg:Ding}--\ref{alg:Mirzal}.  Therefore, we must find effective methods for solving these sub-problems.
Similarly to Lin, we  apply the PG method to solve the sub-problems \eqref{subprobG} -- \eqref{subprobH}.
Algorithm \ref{alg:PG2} contains the main steps of  the PG method for solving the latter and can be straightforwardly
adapted for the former. 

For the sake of simplicity, we denote by $F_H$ the function that we optimize in \eqref{subprobG}, which is actually a simplified version (pure $H$ terms removed) of the objective function from \eqref{ONMF-M} for  $H$  fixed:

$$F_H(G):= \frac12\norm{R-GH}_F^2+\frac{\beta}{2}\norm{G^TG-I}_F^2.
$$ Similarly, for  $G$ is fixed, the objective function from  \eqref{subprobH} will be  denoted by:
$$F_G(H):= \frac12\norm{R-GH}_F^2+\frac{\alpha}{2}\norm{HH^T-I}_F^T.$$

In Algorithm \ref{alg:PG2}, $P$ is the projection operator which projects the new point (matrix) on the cone of non-negative matrices (we simply put negative entries to 0). 
Inequality \eqref{ArmigoRule} shows the Armijo rule to find a suitable step-size guaranteeing a sufficient decrease. Searching for $\lambda_k$ is a time-consuming operation, therefore we strive to do  only a small number of trials for new $\lambda$  in Step 3.1.
Similarly to  Lin  \cite{Lin2007}, we allow for $\lambda$ any positive value. More precisely, we start with $\lambda=1$ and if the Armijo rule \eqref{ArmigoRule} is satisfied, we increase the value of $\lambda$ by dividing it with $\gamma<1$. We repeat this until \eqref{ArmigoRule} is no longer satisfied or the same matrix $H_\lambda$ as in the previous iteration is obtained. 
If the starting $\lambda=1$ does not yield $H_\lambda$ which would satisfy the Armijo rule \eqref{ArmigoRule}, then we decrease it by a factor  $\gamma$ and repeat this until \eqref{ArmigoRule} is satisfied. 
The numerical results obtained using different values of  parameters $\gamma$
(updating  factor  for $\lambda$) and $\sigma$ (parameter to check \eqref{ArmigoRule}) are reported in the following subsections.


\begin{figure} [!htb]\ralg{alg:PG2}
\algrule[1pt]
  {\bf{Algorithm 4. PG method using Armijo rule to solve sub-problem \eqref{subprobH} }} \vspace{-3mm}\\
\algrule[1pt]
\textbf{INPUT:} $ 0<\sigma<1,\gamma<1$, and initial $H^0$.
\begin{enumerate}
\item  Set $k=0$
\item  \textbf{Repeat}
\begin{itemize}
\item[] Find a $\lambda$ (using updating factor $\gamma$) such that for $H_\lambda:=P[H^{k}-\lambda\nabla F_G(H^{k})]$ we have
 \bea  \label{ArmigoRule} F_G(H_\lambda)-F_G(H^{k})\leq\sigma \nabla F_G(H^{k})(H_\lambda-H^{k});\eea
\item[]  Set $H^{k+1}:=H_\lambda$
\item[] Set $k=k+1$;
\end{itemize}
\item \textbf{Until} some stopping  criteria is satisfied.
\end{enumerate}
\textbf{OUTPUT:} $H=H^{k+1}.$
\algrule[1pt]
\vspace{-3mm}
\end{figure}

\subsection{Stopping criteria for Algorithms \ref{alg:PG1} and \ref{alg:PG2}}

As practiced in the literature (e.g. see \cite{Lin-More}), in a constrained optimization problem with the non-negativity constraint on the variable $x$,  
a common condition to check whether a point $x^k$ is close to a stationary point is
\bea  \label{main stop criteri}
\norm{\nabla^P f(x^k)}\leq \varepsilon \norm{\nabla f(x^0)},
\eea
where $f$ is the differentiable function that we try to optimize and $\nabla^P f(x^k)$ is the projected gradient defined as
\bea \label{def of projected gradient}
\nabla^P f(x)_i=\begin{cases} 
\nabla f(x)_i,~~~~~~~~~~~~~~~~~\text{if}~ x_{i}>0,\\
\min\{0,\nabla f(x)_i\},~~~~~\text{if} ~  x_{i}=0,
\end{cases}
\eea
and $\varepsilon$ is a small positive tolerance.
For Algorithm \ref{alg:PG1}, \eqref{main stop criteri} becomes
\bea  \label{stoping criteri for ONMF}
\norm{\nabla^P F\br{G^k,H^k}}_F\leq \norm{\varepsilon {\nabla F\br{G^0,H^0}}}_F.
\eea
 We impose a time limit in seconds and a maximum
 number of iterations for Algorithm \ref{alg:PG2} as well.
 Following \cite{Lin2007}, we also define stopping conditions for the sub-problems. The
 matrices $G^{ k+1}$ and $H^{k+1}$ returned by Algorithm \ref{alg:PG2}, respectively, must satisfy
\bea \label{stoping criteria for subprobs}
\begin{array}{ccc}
\norm{\nabla_{G}^PF\br{G^{k+1}, H^k}}_F\leq \bar{\varepsilon}_G,\\
\norm{\nabla_{H}^PF\br{G^{k+1}, H^{k+1}}}_F\leq \bar{\varepsilon}_H,
 \end{array}
\eea
where
\bea \label{def epsG and epsH}
 \bar{\varepsilon}_G= \bar{\varepsilon}_H=\max\{10^{-7}, \varepsilon\}\norm{\nabla F\br{G^{0}, H^0}}_F,
\eea
and $\varepsilon$ is the same tolerance used in \eqref{stoping criteri for ONMF}. If the PG method for solving the sub-problem \eqref{subprobG} or \eqref{subprobH} stops
 after the first iteration, then we decrease the stopping tolerance as follows:
\beq  \label{corrected tolerance for subprobs}
\bar{\varepsilon}_G\longleftarrow \tau\bar{\varepsilon}_G,~~~~ \bar{\varepsilon}_H\longleftarrow \tau\bar{\varepsilon}_H,
\eeq
where $\tau$ is a constant smaller then 1.
\section{Numerical results}
\label{sec:num}
In this section we demonstrate, how the PG method described in Sect. \ref{sec:PGM}, performs compared to the MU-based algorithms of Ding and Mirzal, which were described in Subsections \ref{subsec:DING} and \ref{subsec:Mirzal}, respectively.

\subsection{Artificial data}
We created two sets of synthetic data using MATLAB \cite{matlab}.
The first set we call bi-orthonormal set (BION).
It consists of instances of matrix $R\in\R_+^{n\times n}$, which were created as products of $G$ and $H$, where $G\in\R_+^{n\times k}$
has orthonormal columns while  $H \in \R_+^{k\times n}$ has orthonormal rows.
We created five instances of $R$, for each pair $(n,k_1)$ and $(n,k_2)$ from Table \ref{tab:n_k}.

Matrices $G$ were created in two phases: firstly, we randomly (uniform distribution) selected a position in each row; 
secondly, we selected a random number from $(0,1)$ (uniform distribution) for the selected position in each row.
Finally, if it happens that after this procedure some column of $G$ is zero or has a norm below $10^{-8}$, we find the first non-zero element in the largest  column of $G$ (according to Euclidean norm) and move it into the zero column. We created $H$ similarly.
\begin{table}
	$$
	\begin{array}{c||c|c|c|c|c|c}
		n & 50 & 100  & 200 & 500 & 1000\\\hline
		k_1 & 10 &20 & 40 & 100 &200\\\hline
		k_2 & 20 &40 & 80 & 200 & 400
	\end{array}$$
	 \caption{Paris $(n,k)$ for which we created UNION and BION  datasets}\label{tab:n_k}
\end{table}
Each triple $(R,G,H)$ was saved as a triple of {\tt txt} files. For example, {\tt NMF\_BIOG\_data\_R\_n=200\_k=80\_id=5.txt} contains $200\times 200$ matrix $R$ obtained by 
multiplying matrices $G\in\R^{200\times 80}$ and $H\in \R^{80\times 200}$, which were generated as explained above. With {\tt id=5}, we denote that this is a 5th matrix corresponding to this pair $(n,k)$.
The second set contains similar data to BION, but only one factor ($G$) is orthonormal, while the other ($H$) is non-negative but not necessarily orthonormal.
We call this dataset uni-orthonormal (UNION).
All computations are done using MATLAB \cite{matlab} and a high performance computer available at Faculty of Mechanical Engineering of University of Ljubljana.
This is Intel Xeon X5670 (1536 hyper-cores) HPC cluster and an E5-2680 V3 (1008 hyper-cores) DP cluster, with an
IB QDR interconnection, 164 TB of LUSTRE storage, 4.6 TB RAM and with 24 TFlop/s performance. 

\subsection{Numerical results for UNION}

In this subsection, we present numerical results, obtained by Ding's, Mirzal's, and our algorithm for a uni-orthogonal problem \eqref{ONMF}, using the UNION data, introduced in the previous subsection.
We have adapted the last two algorithms (Algorithms \ref{alg:Mirzal}, \ref{alg:PG1})
for UNION data by setting $\alpha=0$ in the problem formulation \eqref{biONMF} and in all formulas underlying these two algorithms. 

Recall that for UNION data we have 
for each pair $n,k$ from Table \ref{tab:n_k} five symmetric matrices $R$ for which we try to solve \eqref{ONMF} by Algorithms \ref{alg:Ding}, \ref{alg:Mirzal} and 
\ref{alg:PG1}. Note that all these algorithms demand as input the internal dimension $k$, i.e. the number of columns of factor $G$, which is in general not known in advance. Even though, we know this dimension by construction for UNION data, we tested the algorithms using internal dimensions $p$ equal to $20\%,~40\%,\ldots,100\%$ of $k$.
For $p=k$, we know the optimum of the problem, which is 0, so for this case we can also estimate how good are the tested algorithms in terms of finding the global optimum.

The first question we had to answer was which value of $\beta$ to use in Mirzal's and PG algorithms.
It is obvious that larger values of $\beta$ moves the focus from optimizing the $\text{RSE}$ to guaranteeing the orthonormality, i.e., feasibility for the original problem.
We decided not to fix the value of $\beta$ but to run both algorithms for $\beta\in\{1,10,100,1000\}$ and report the results.

For each solution pair $G,H$ returned by all algorithms, the non-negativity constraints are held by the construction of algorithms, so we only need to consider deviation of $G$ from orthonormality, which we call \emph{infeasibility} and define it as
\bea \label{def:vio}
\mbox{infeas}_G:= \frac{\norm{G^TG-I}_F}{1+\norm{I}_F}.
\eea

The  computational results that follow in the rest of this subsection were obtained by setting the tolerance in the stopping criterion to $\varepsilon=10^{-10}$, the maximum number of iterations to $1000$ in Algorithm \ref{alg:PG1} and to 20 in Algorithm \ref{alg:PG2}.
We also set  a time limit to  $3600$ seconds. Additionally,  for  $\sigma$ and $\gamma$ (updating parameter for $\lambda$ in Algorithm \ref{alg:PG2}) 
we choose  $0.001$ and $0.1$, respectively. Finally, for $\tau$ from \eqref{corrected tolerance for subprobs} we set a value of $0.1$.

In general, Algorithm \ref{alg:PG1} converges to a solution in early iterations and the norm of the projected gradient falls below the tolerance shortly after running the algorithm.

\begin{figure}[h!]
\centering
  \begin{subfigure}{6cm} 
    \centering
    \includegraphics[width=5cm]{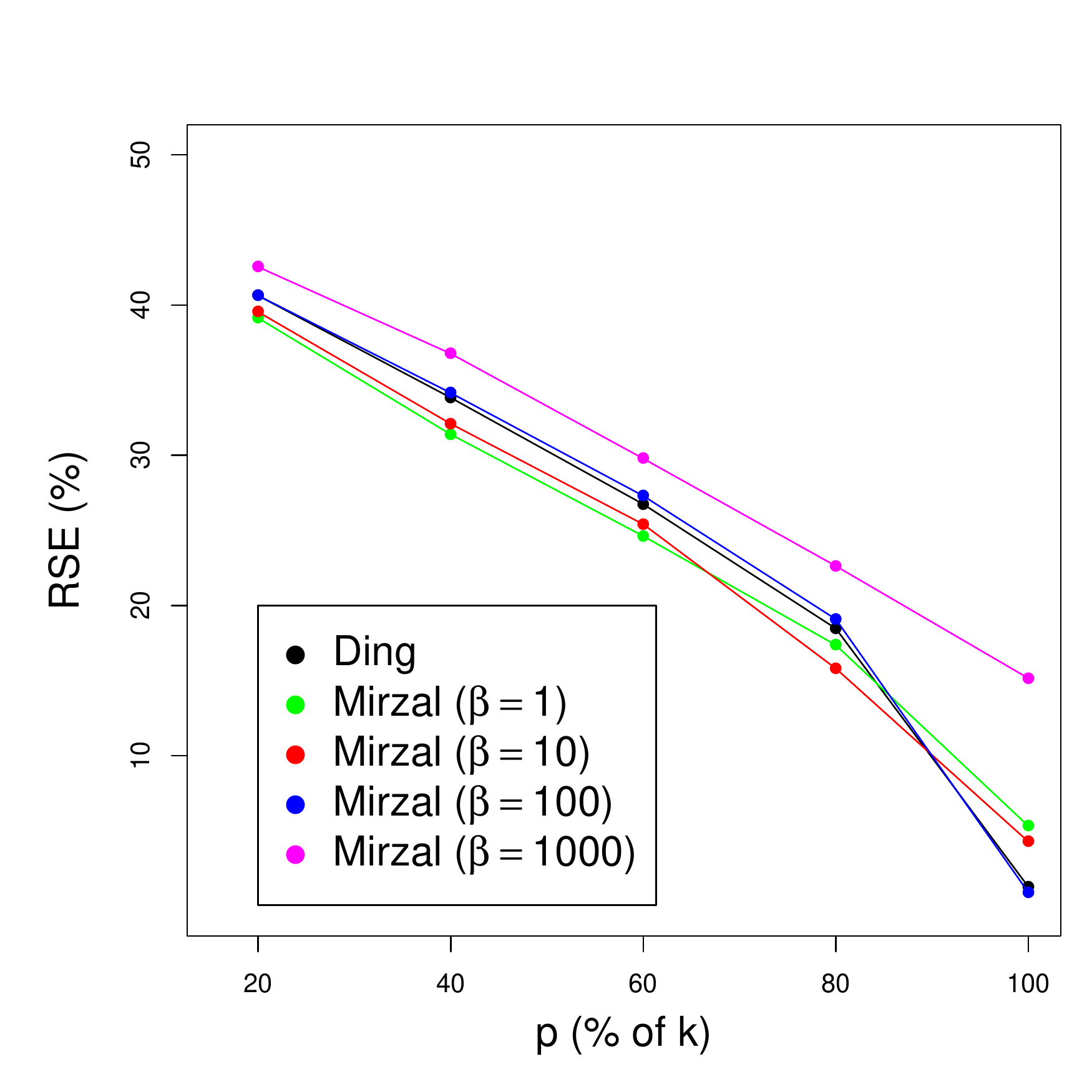}
    \caption{Values of $\text{RSE}$ for different values of $\beta$ obtained by  Algorithms \ref{alg:Ding} and \ref{alg:Mirzal} for $n=100$}\label{fig:rse_infeas_union_c}
  \end{subfigure}
  \begin{subfigure}{0.7cm}
  ~
  \end{subfigure}
  \begin{subfigure}{6cm}
    \centering\includegraphics[width=5cm]{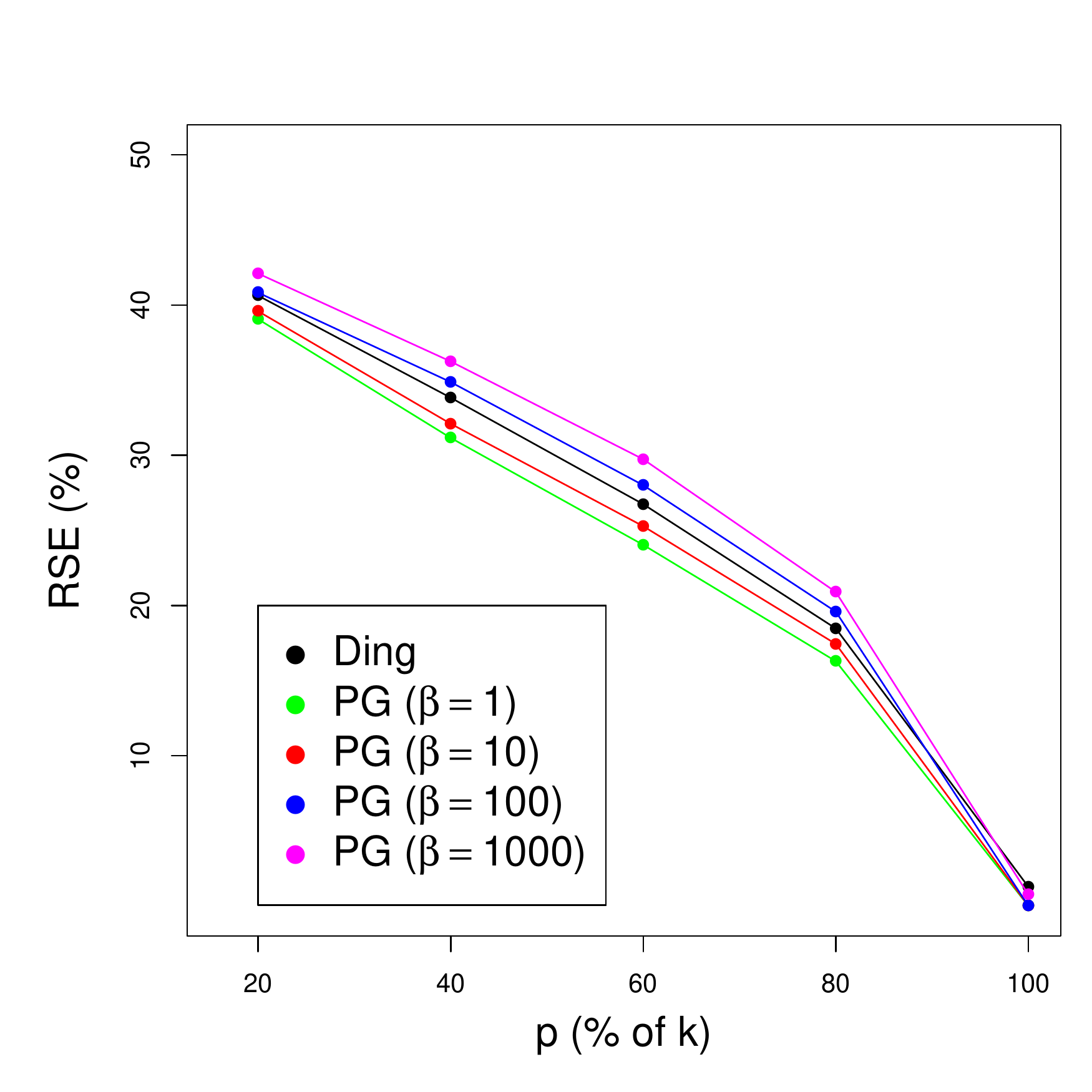}
    \caption{Values of $\text{RSE}$ for different values of $\beta$ obtained by  Algorithms \ref{alg:Ding} and  \ref{alg:PG1} for $n=100$}\label{fig:rse_infeas_union_d}  
    \end{subfigure}

  \begin{subfigure}{6cm} 
    \centering
    \includegraphics[width=5cm]{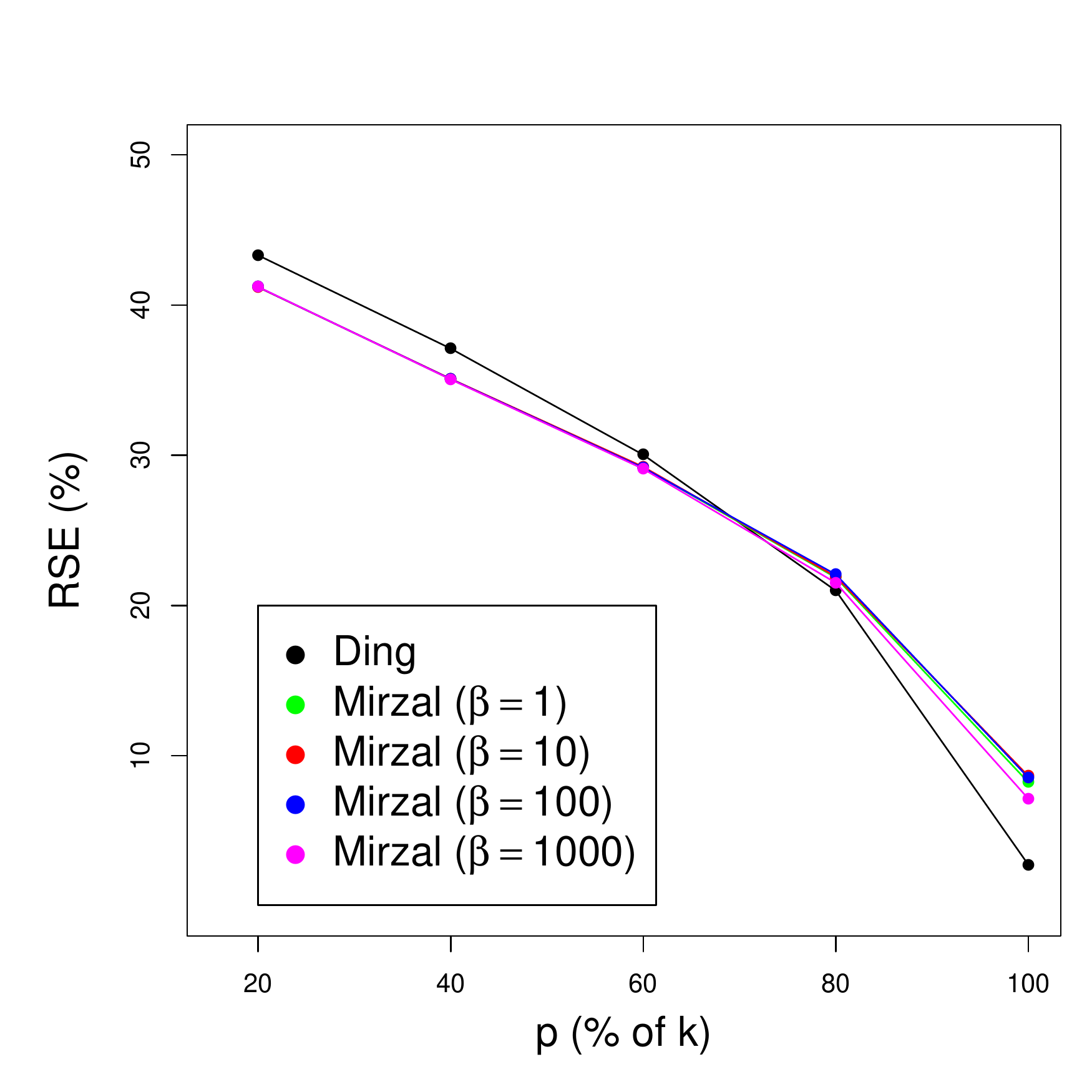}
    \caption{Values of $\text{RSE}$ for different values of $\beta$ obtained by  Algorithms \ref{alg:Ding} and  \ref{alg:Mirzal} for $n=500$}\label{fig:rse_infeas_union_g}
  \end{subfigure}
  \begin{subfigure}{0.7cm}
  ~
  \end{subfigure}
  \begin{subfigure}{6cm}
    \centering\includegraphics[width=5cm]{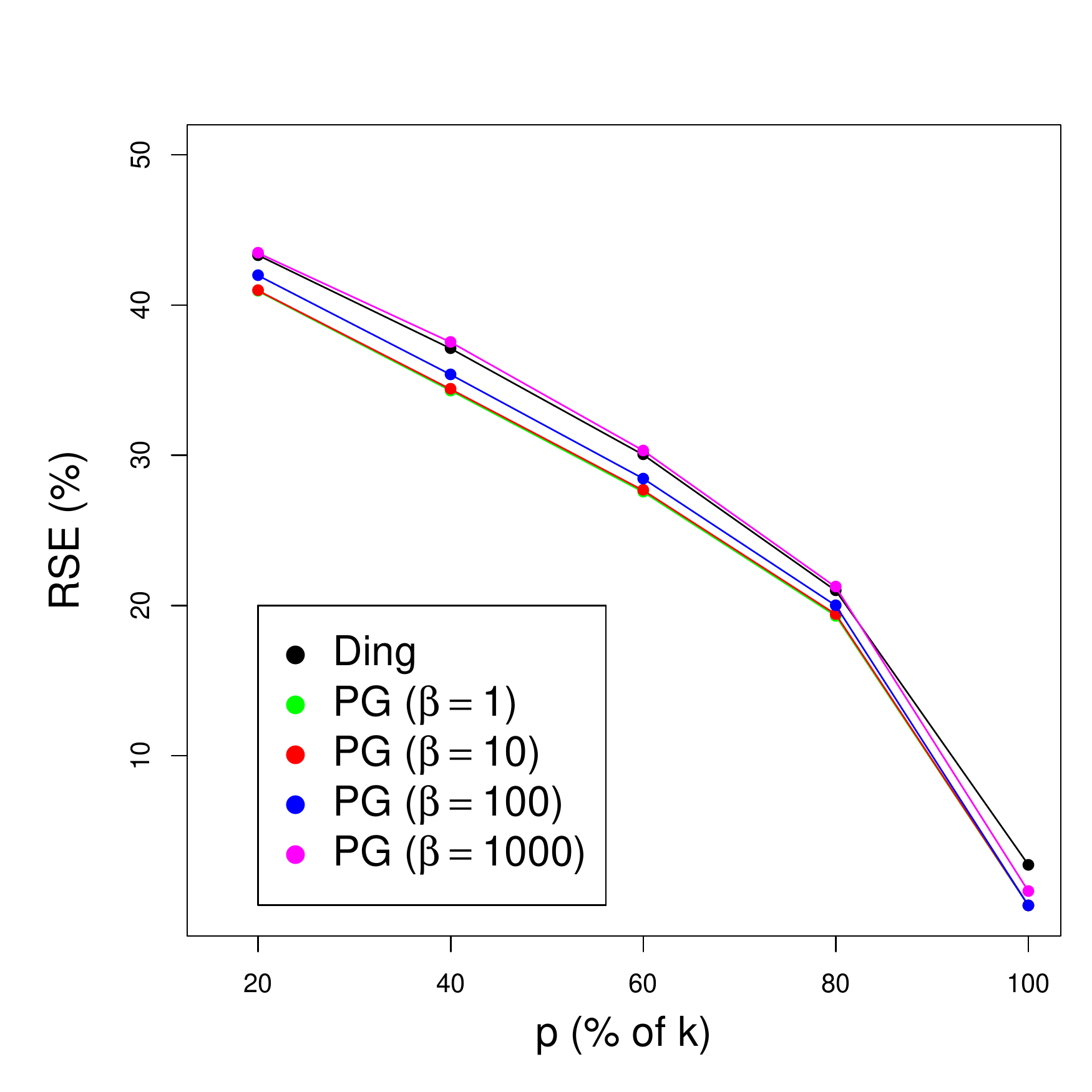}
    \caption{Values of $\text{RSE}$ for different values of $\beta$ obtained by  Algorithms \ref{alg:Ding} and  \ref{alg:PG1} for $n=500$}\label{fig:rse_infeas_union_h}  
    \end{subfigure}

  \begin{subfigure}{6cm} 
    \centering
    \includegraphics[width=5cm]{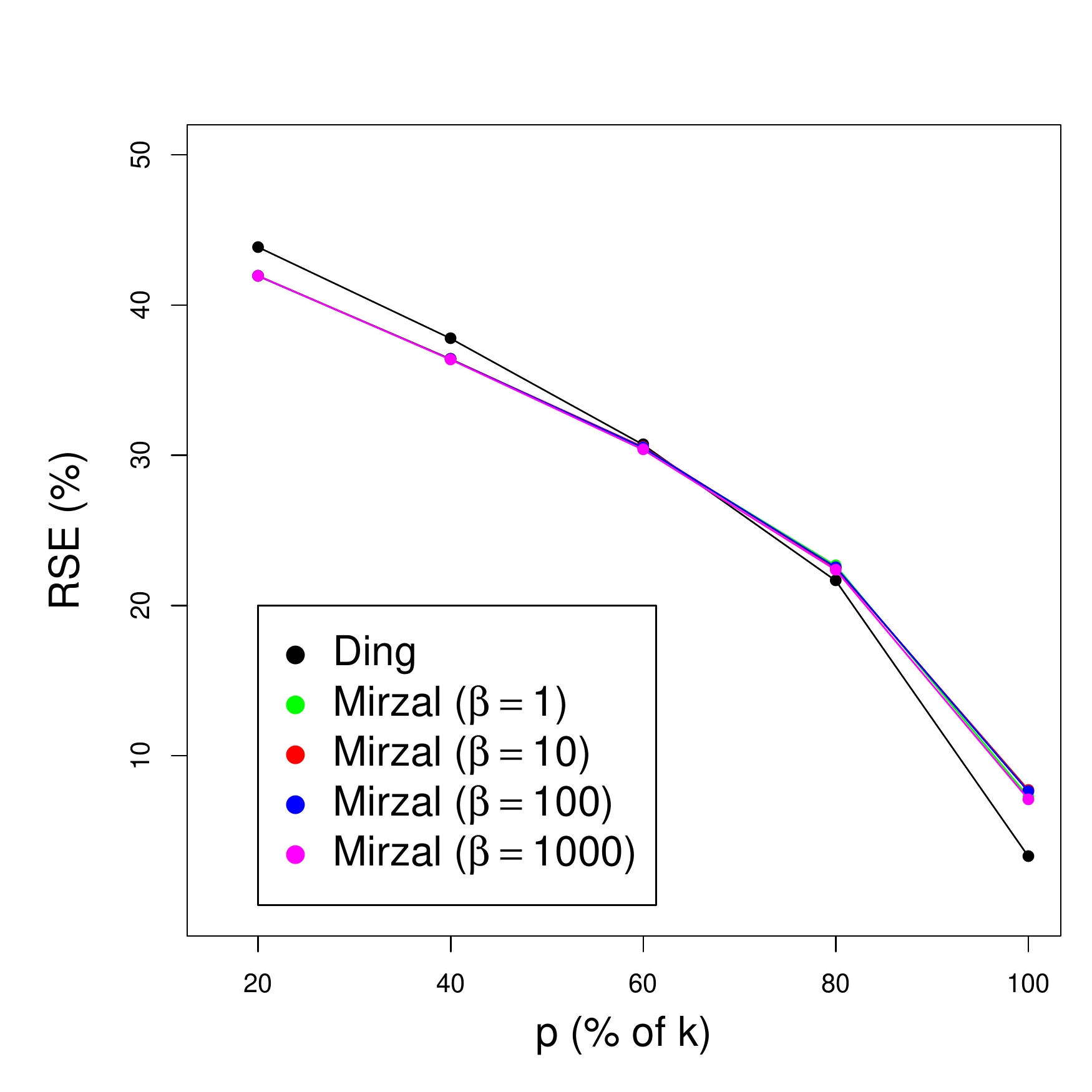}
    \caption{Values of $\text{RSE}$ for different values of $\beta$ obtained by  Algorithms \ref{alg:Ding} and  \ref{alg:Mirzal} for $n=1000$}\label{fig:rse_infeas_union_i}
  \end{subfigure}
    \begin{subfigure}{0.7cm}
  ~
  \end{subfigure}
  \begin{subfigure}{6cm}
    \centering\includegraphics[width=5cm]{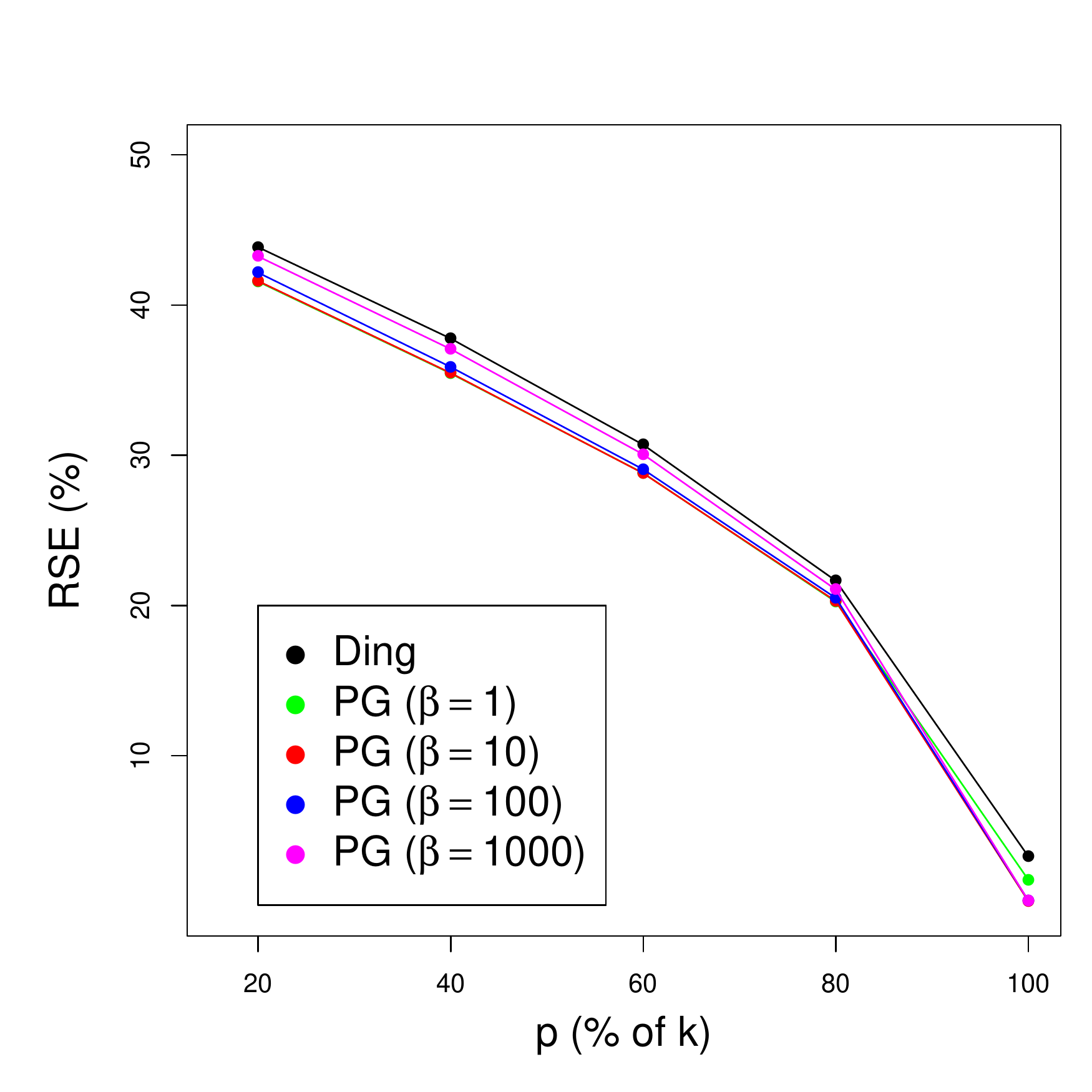}
    \caption{Values of $\text{RSE}$ for different values of $\beta$ obtained by  Algorithms \ref{alg:Ding} and  \ref{alg:PG1} for $n=1000$}\label{fig:rse_infeas_union_j}  
    \end{subfigure}

    \caption{This figure depicts data from Table \ref{tab:RSE_vs_beta}. It contains six plots which illustrate the quality of Algorithms  \ref{alg:Ding}, \ref{alg:Mirzal} and \ref{alg:PG1} regarding $\text{RSE}$ on UNION instances with  $n=100,500,1000$, for $\beta\in\{1,10,100,1000\}$. We can see that regarding $\text{RSE}$ the performance of these algorithms on this dataset does not differ a lot. As expected, larger values of $\beta$ yield larger values of $\text{RSE}$, but the differences are rather small.  However, when $p$ approached 100 \% of $k$, Algorithm  \ref{alg:PG1} comes closest to the global optimum  $\text{RSE}=0$.}\label{fig:rse_vs_beta_a}
\end{figure}

\begin{figure}[h!]
\centering
  \begin{subfigure}{6cm} 
    \centering
    \includegraphics[width=5cm]{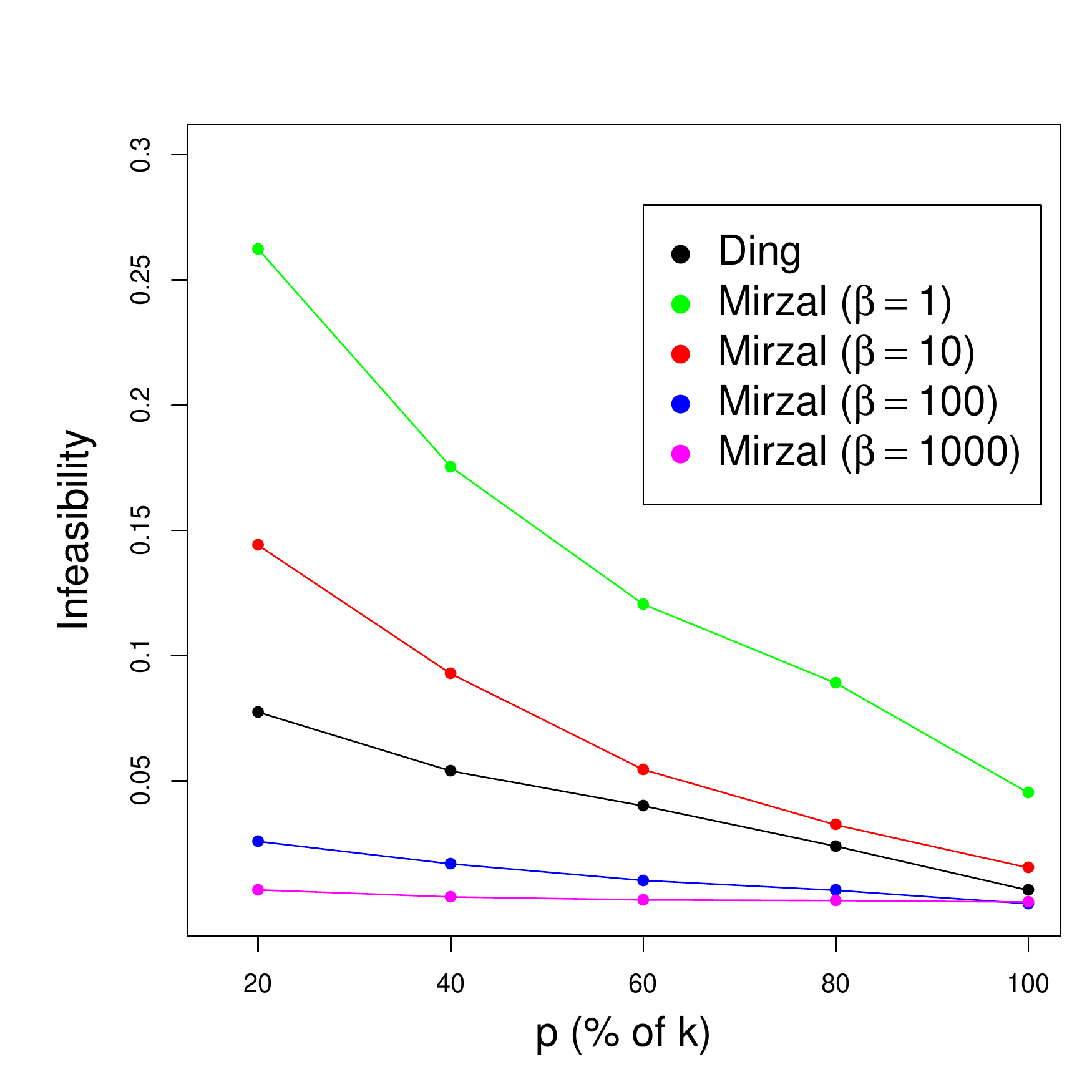}
    \caption{Values of $\mbox{infeas}_G$ for different values of $\beta$ obtained by  Algorithms \ref{alg:Ding} and \ref{alg:Mirzal} for $n=100$}\label{fig:infeas_vs_beta_union_a}
  \end{subfigure}
  \begin{subfigure}{0.7cm}
  ~
  \end{subfigure}
  \begin{subfigure}{6cm}
    \centering\includegraphics[width=5cm]{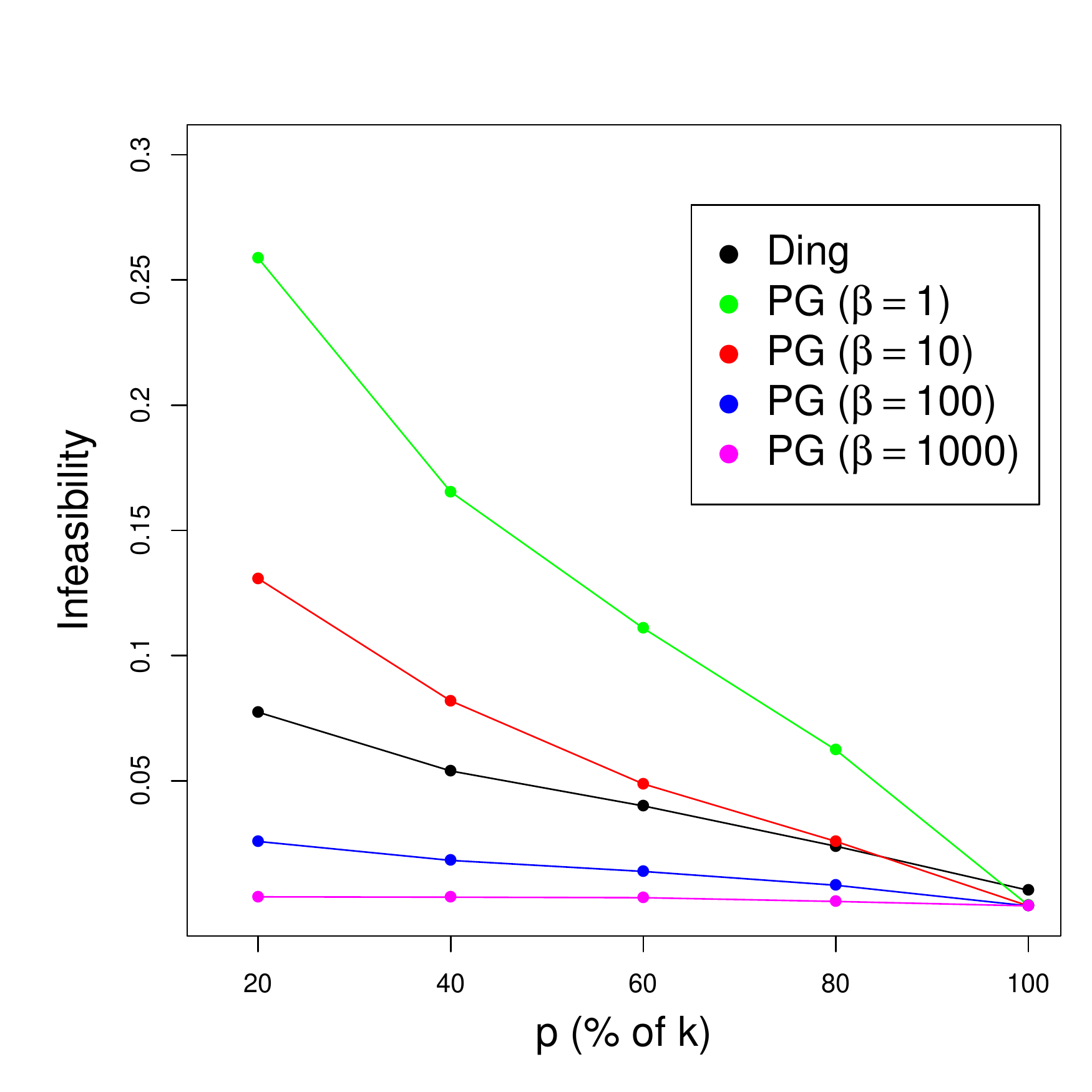}
    \caption{Values of $\mbox{infeas}_G$ for different values of $\beta$ obtained by  Algorithms \ref{alg:Ding} and  \ref{alg:PG1} for $n=100$}\label{fig:infeas_vs_beta_union_b}  
    \end{subfigure}

  \begin{subfigure}{6cm} 
    \centering
    \includegraphics[width=5cm]{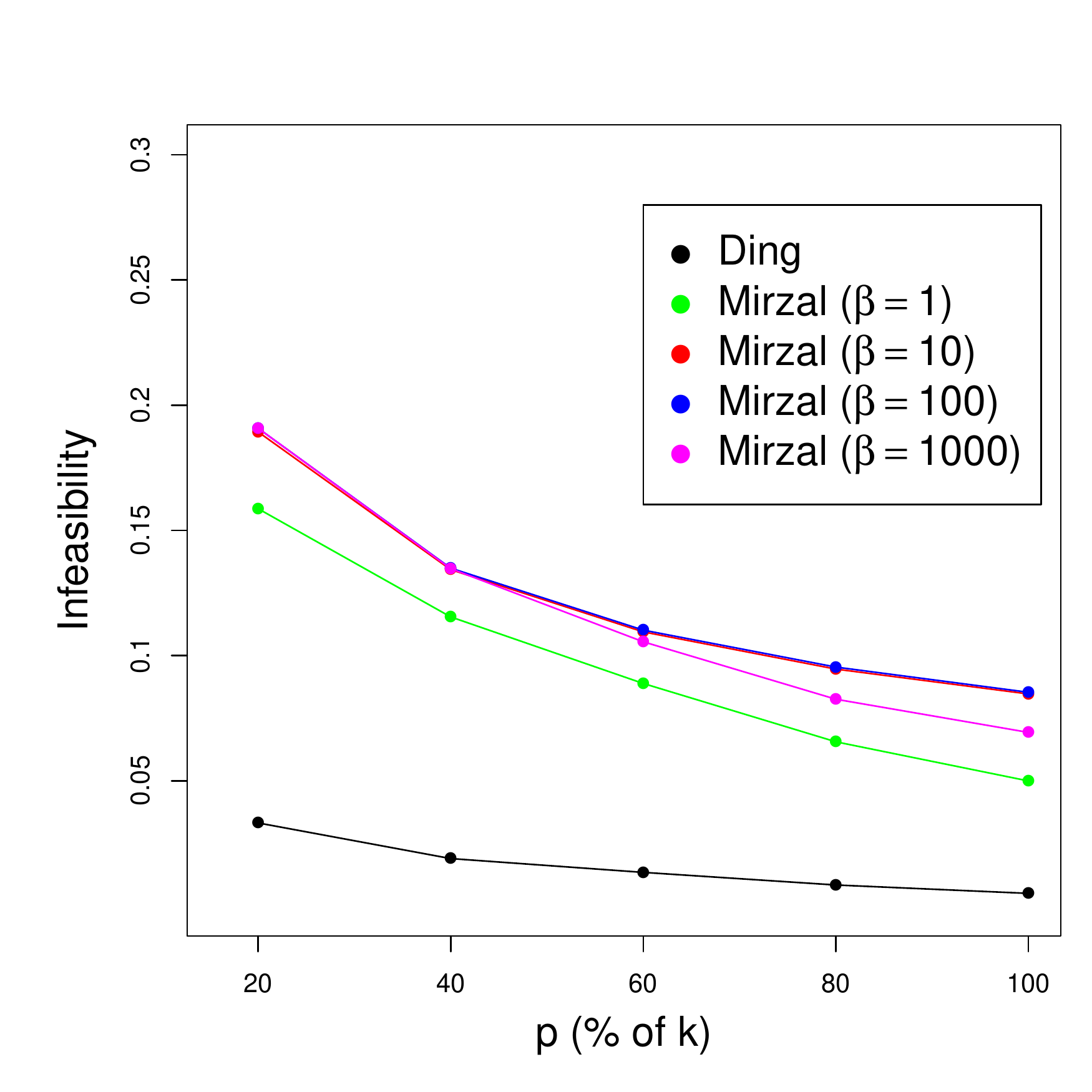}
    \caption{Values of $\mbox{infeas}_G$ for different values of $\beta$ obtained by  Algorithms \ref{alg:Ding} and  \ref{alg:Mirzal} for $n=500$}\label{fig:infeas_vs_beta_union_c}
  \end{subfigure}
  \begin{subfigure}{0.7cm}
  ~
  \end{subfigure}
  \begin{subfigure}{6cm}
    \centering\includegraphics[width=5cm]{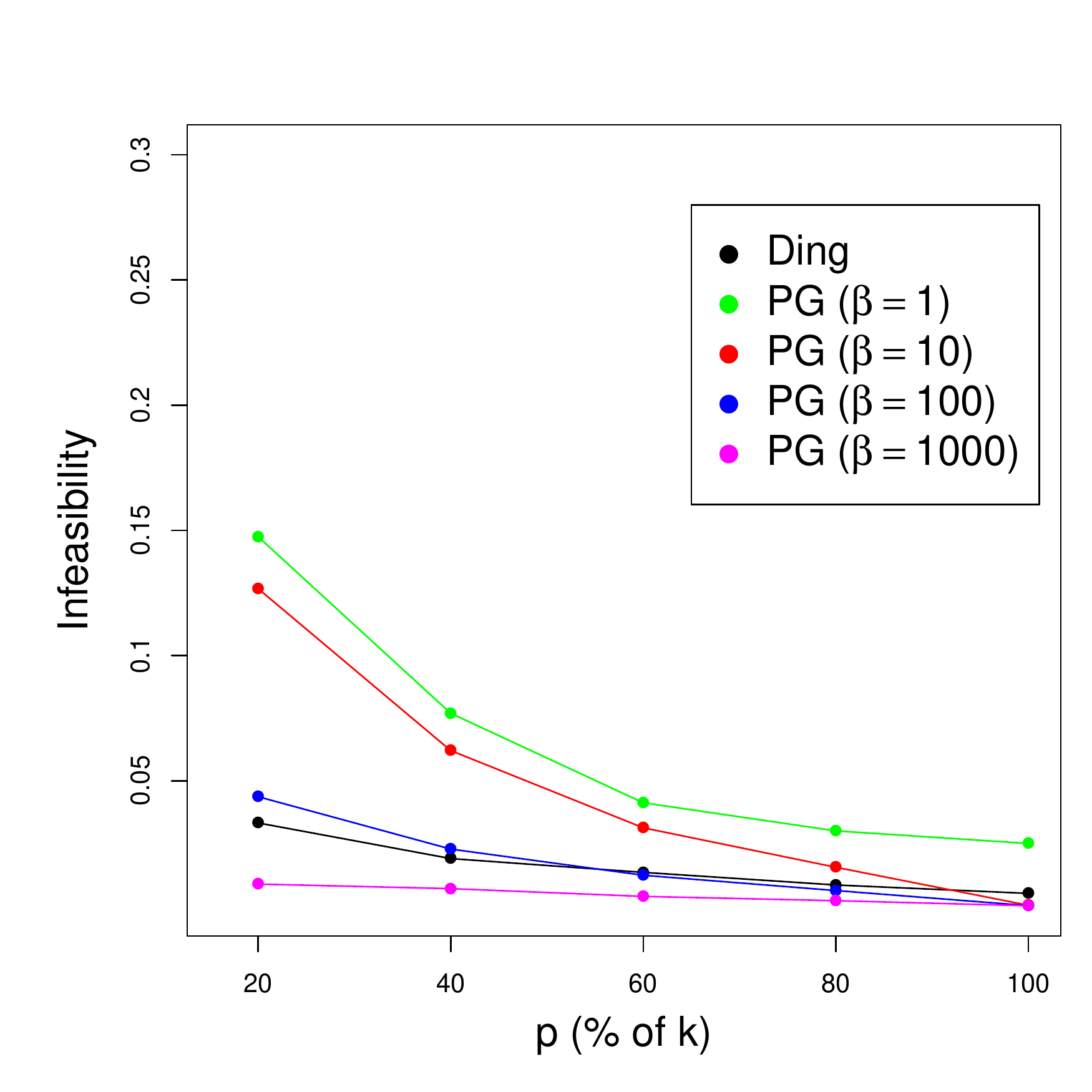}
    \caption{Values of $\mbox{infeas}_G$ for different values of $\beta$ obtained by  Algorithms \ref{alg:Ding} and  \ref{alg:PG1} for $n=500$}\label{fig:infeas_vs_beta_union_d}  
    \end{subfigure}

  \begin{subfigure}{6cm} 
    \centering
    \includegraphics[width=5cm]{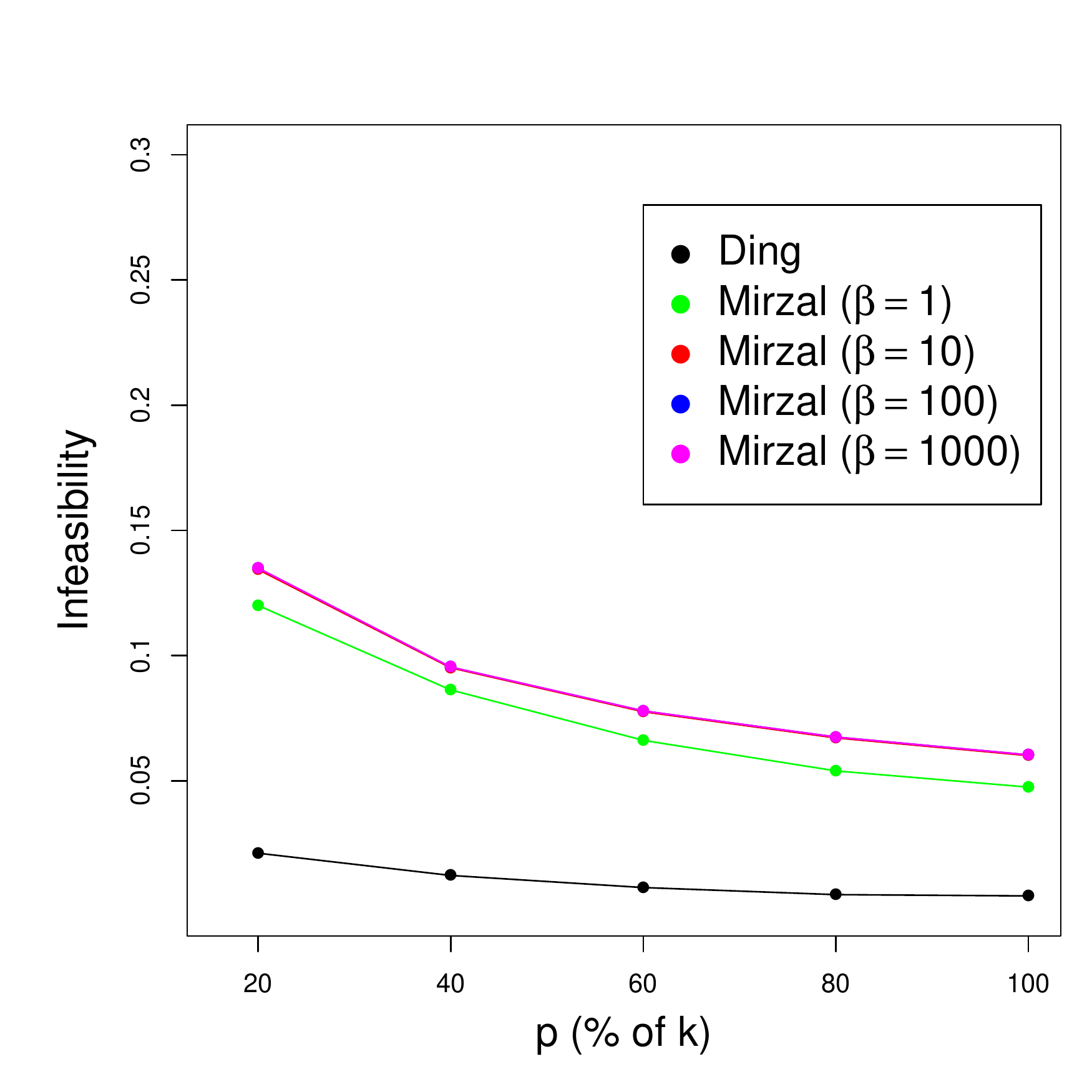}
    \caption{Values of $\mbox{infeas}_G$ for different values of $\beta$ obtained by  Algorithms \ref{alg:Ding} and  \ref{alg:Mirzal} for $n=1000$}\label{fig:infeas_vs_beta_union_e}
  \end{subfigure}
    \begin{subfigure}{0.7cm}
  ~
  \end{subfigure}
  \begin{subfigure}{6cm}
    \centering\includegraphics[width=5cm]{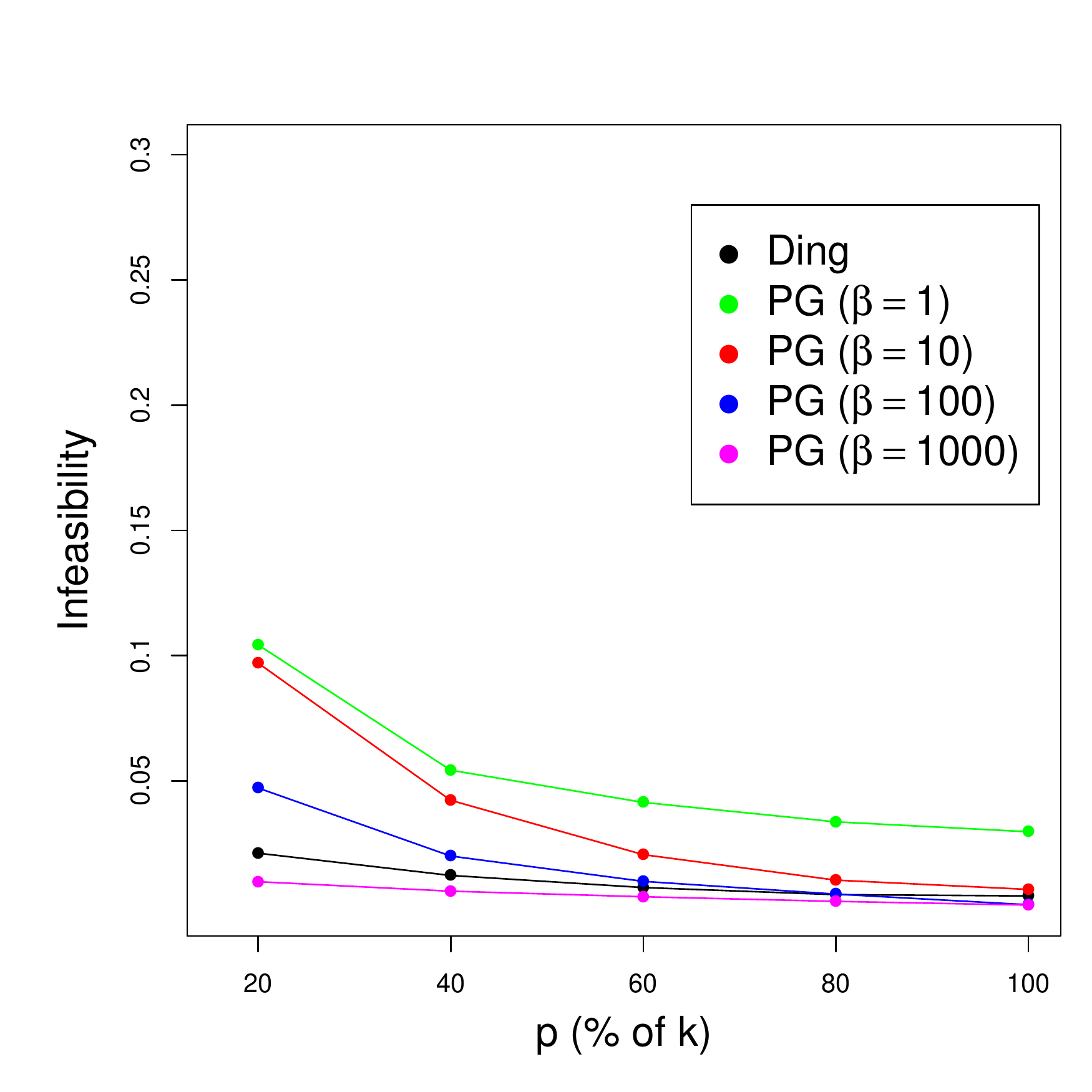}
    \caption{Values of $\mbox{infeas}_G$ for different values of $\beta$ obtained by  Algorithms \ref{alg:Ding} and  \ref{alg:PG1} for $n=1000$}\label{fig:infeas_vs_beta_union_f}  
    \end{subfigure}

    \caption{This figure depicts data from Table \ref{tab:infeas_vs_beta}. It contains six plots which illustrate the quality of Algorithms \ref{alg:Ding}, \ref{alg:Mirzal} and \ref{alg:PG1} regarding infeasibility  on UNION instances with  $n=100,500,1000$, for $\beta\in\{1,10,100,1000\}$. We can see that regarding infeasibility the performance of these algorithms on this dataset does not differ a lot. As expected, larger values of $\beta$ yield smaller values of $\mbox{infeas}_G$, but the differences are rather small.}\label{fig:infeas_vs_beta}
\end{figure}

Results in Tables \ref{tab:RSE_vs_beta} and \ref{tab:infeas_vs_beta} and their visualisations on Figures \ref{fig:rse_infeas_union_c}--Figures \ref{fig:rse_infeas_union_j} and on 
Figures \ref{fig:infeas_vs_beta_union_a}--Figures \ref{fig:infeas_vs_beta_union_f} confirm expectations. 
More precisely, we can see that the smaller the value of $\beta$, the better $\text{RSE}$. Likewise, the larger the value of $\beta$, the smaller the infeasibility 
$\mbox{infeas}_G$.
In practice, we want to reach both criteria: small $\text{RSE}$ and small infeasibility, so some compromise should be made. If $\text{RSE}$ is more important than infeasibility, we choose the smaller value of $\beta$ and vice versa.
We can also observe that regarding $\text{RSE}$ the three compared algorithms do not differ a lot. However, when the input dimension $p$ approaches the real inner dimension $k$,  Algorithm \ref{alg:PG1} comes closest to the global optimum $\text{RSE}=0$.
The situation with infeasibility is a bit different. While Algorithm \ref{alg:Ding} performs very well in all instances, Algorithm \ref{alg:Mirzal} reaches better feasibility for smaller values of $n$. Algorithm \ref{alg:PG1}  outperforms the others for $\beta =1000$.

\begin{table}[h!]
\centering \scriptsize
\begin{tabular}{r|r||r|rrrr|rrrr} 
 \multicolumn{1}{c|}{\multirow{2}{*}{$n$}} & \multicolumn{1}{c||}{$p$} & $\text{RSE}$ of  & \multicolumn{4}{c|}{$\text{RSE}$ of Alg.  \ref{alg:Mirzal}} & \multicolumn{4}{c}{$\text{RSE}$ of Alg.  \ref{alg:PG1}}\\
& \multicolumn{1}{c||}{(\% of $k$)} & Alg.  \ref{alg:Ding} &$\beta=1$ & $\beta=10$& $\beta=100$& $\beta=1000$& $\beta=1$& $\beta=10$& $\beta=100$& $\beta=1000$
\\ 
  \hline
  50 & 40 & 0.3143 & 0.2965 & 0.3070 & 0.3329 & 0.3898 & 0.2963 & 0.3081 & 0.3425 & 0.3508 \\ 
  50 & 60 & 0.2348 & 0.2227 & 0.2356 & 0.2676 & 0.3459 & 0.2201 & 0.2382 & 0.2733 & 0.2765 \\ 
  50 & 80 & 0.1738 & 0.1492 & 0.1634 & 0.1894 & 0.3277 & 0.1468 & 0.1620 & 0.1953 & 0.2053 \\ 
  50 & 100 & 0.0002 & 0.0133 & 0.0004 & 0.0932 & 0.2973 & 0.0000 & 0.0000 & 0.0000 & 0.0000 \\ 
  100 & 20 & 0.4063 & 0.3914 & 0.3955 & 0.4063 & 0.4254 & 0.3906 & 0.3959 & 0.4083 & 0.4210 \\ 
  100 & 40 & 0.3384 & 0.3139 & 0.3210 & 0.3415 & 0.3677 & 0.3116 & 0.3210 & 0.3488 & 0.3625 \\ 
  100 & 60 & 0.2674 & 0.2462 & 0.2541 & 0.2730 & 0.2978 & 0.2403 & 0.2528 & 0.2801 & 0.2974 \\ 
  100 & 80 & 0.1847 & 0.1737 & 0.1581 & 0.1909 & 0.2263 & 0.1629 & 0.1744 & 0.1959 & 0.2090 \\ 
  100 & 100 & 0.0126 & 0.0532 & 0.0427 & 0.0089 & 0.1515 & 0.0000 & 0.0000 & 0.0000 & 0.0075 \\ 
  200 & 20 & 0.4213 & 0.4024 & 0.4077 & 0.4080 & 0.4257 & 0.4005 & 0.4032 & 0.4162 & 0.4337 \\ 
  200 & 40 & 0.3562 & 0.3315 & 0.3398 & 0.3401 & 0.3647 & 0.3270 & 0.3313 & 0.3497 & 0.3738 \\ 
  200 & 60 & 0.2845 & 0.2675 & 0.2746 & 0.2748 & 0.2955 & 0.2573 & 0.2617 & 0.2812 & 0.3061 \\ 
  200 & 80 & 0.1959 & 0.1958 & 0.2013 & 0.1996 & 0.2085 & 0.1773 & 0.1819 & 0.1960 & 0.2133 \\ 
  200 & 100 & 0.0191 & 0.0753 & 0.0632 & 0.0622 & 0.0415 & 0.0000 & 0.0000 & 0.0069 & 0.0181 \\ 
  500 & 20 & 0.4332 & 0.4120 & 0.4119 & 0.4120 & 0.4121 & 0.4092 & 0.4096 & 0.4197 & 0.4346 \\ 
  500 & 40 & 0.3711 & 0.3506 & 0.3509 & 0.3507 & 0.3505 & 0.3430 & 0.3440 & 0.3537 & 0.3753 \\ 
  500 & 60 & 0.3003 & 0.2919 & 0.2923 & 0.2916 & 0.2909 & 0.2756 & 0.2766 & 0.2845 & 0.3031 \\ 
  500 & 80 & 0.2098 & 0.2186 & 0.2192 & 0.2207 & 0.2151 & 0.1931 & 0.1941 & 0.1999 & 0.2122 \\ 
  500 & 100 & 0.0273 & 0.0822 & 0.0864 & 0.0853 & 0.0713 & 0.0002 & 0.0003 & 0.0002 & 0.0097 \\ 
  1000 & 20 & 0.4386 & 0.4195 & 0.4194 & 0.4193 & 0.4195 & 0.4156 & 0.4160 & 0.4216 & 0.4324 \\ 
  1000 & 40 & 0.3777 & 0.3641 & 0.3640 & 0.3638 & 0.3637 & 0.3545 & 0.3548 & 0.3588 & 0.3707 \\ 
  1000 & 60 & 0.3070 & 0.3047 & 0.3055 & 0.3051 & 0.3036 & 0.2881 & 0.2880 & 0.2906 & 0.3006 \\ 
  1000 & 80 & 0.2164 & 0.2265 & 0.2248 & 0.2254 & 0.2236 & 0.2024 & 0.2029 & 0.2050 & 0.2106 \\ 
  1000 & 100 & 0.0329 & 0.0725 & 0.0772 & 0.0761 & 0.0709 & 0.0173 & 0.0030 & 0.0035 & 0.0035
\end{tabular}
\caption{In this table we demonstrate how good $\text{RSE}$ is achieved by Algorithms \ref{alg:Ding}, \ref{alg:Mirzal} and \ref{alg:PG1} on UNION dataset. 
For each  $n\in\{50,100,200,500,1000\}$  we take all 10 matrices $R$ (five of them corresponding to $k=0.2n$ and five to $k=0.4n$). We run all three algorithms on these matrices with inner dimensions $p\in \{0.2k,0.4k,\ldots,1.0k\}$
with all possible values of $\beta\in\{1,10,100,1000\}$. 
Each row represents the average (arithmetic mean value) $\text{RSE}$ obtained on instances corresponding 
to given $n$. For example, the last row shows the average value of $\text{RSE}$ in 10 instances of dimension 1000 (five of them corresponding to $k=200$ and five to $k=400$) obtained by all three algorithms for all four values of $\beta$, which were run with the input dimension $p=k$.
} 
\label{tab:RSE_vs_beta}
\end{table}

\begin{table}[h!]
\centering \scriptsize
\begin{tabular}{r|r||r|rrrr|rrrr} 
 \multicolumn{1}{c|}{\multirow{2}{*}{$n$}} & \multicolumn{1}{c||}{$p$} & Infeas. of  & \multicolumn{4}{c|}{Infeas. of Alg.  \ref{alg:Mirzal}} & \multicolumn{4}{c}{Infeas. of Alg.  \ref{alg:PG1}}\\
& \multicolumn{1}{c||}{(\% of $k$)} & Alg.  \ref{alg:Ding} &$\beta=1$ & $\beta=10$& $\beta=100$& $\beta=1000$& $\beta=1$& $\beta=10$& $\beta=100$& $\beta=1000$
\\ 
  \hline
50 & 20 & 0.0964 & 0.2490 & 0.0924 & 0.0155 & 0.0038 & 0.2298 & 0.0909 & 0.0154 & 0.0022 \\ 
  50 & 40 & 0.0740 & 0.1886 & 0.0676 & 0.0131 & 0.0040 & 0.1845 & 0.0670 & 0.0135 & 0.0023 \\ 
  50 & 60 & 0.0553 & 0.1324 & 0.0465 & 0.0068 & 0.0040 & 0.1245 & 0.0440 & 0.0091 & 0.0015 \\ 
  50 & 80 & 0.0324 & 0.0964 & 0.0241 & 0.0053 & 0.0034 & 0.0789 & 0.0250 & 0.0069 & 0.0020 \\ 
  50 & 100 & 0.0023 & 0.0257 & 0.0022 & 0.0023 & 0.0039 & 0.0000 & 0.0000 & 0.0000 & 0.0000 \\ 
  100 & 20 & 0.0774 & 0.2624 & 0.1441 & 0.0258 & 0.0064 & 0.2588 & 0.1308 & 0.0258 & 0.0036 \\ 
  100 & 40 & 0.0539 & 0.1754 & 0.0928 & 0.0168 & 0.0036 & 0.1654 & 0.0819 & 0.0182 & 0.0035 \\ 
  100 & 60 & 0.0400 & 0.1205 & 0.0545 & 0.0102 & 0.0024 & 0.1109 & 0.0487 & 0.0138 & 0.0033 \\ 
  100 & 80 & 0.0239 & 0.0890 & 0.0324 & 0.0062 & 0.0022 & 0.0623 & 0.0258 & 0.0083 & 0.0018 \\ 
  100 & 100 & 0.0062 & 0.0452 & 0.0153 & 0.0009 & 0.0016 & 0.0002 & 0.0000 & 0.0000 & 0.0000 \\ 
  200 & 20 & 0.0584 & 0.2157 & 0.1437 & 0.1433 & 0.0054 & 0.2087 & 0.1512 & 0.0348 & 0.0074 \\ 
  200 & 40 & 0.0356 & 0.1379 & 0.1004 & 0.1000 & 0.0036 & 0.1240 & 0.0806 & 0.0207 & 0.0053 \\ 
  200 & 60 & 0.0260 & 0.0955 & 0.0791 & 0.0793 & 0.0031 & 0.0754 & 0.0434 & 0.0143 & 0.0047 \\ 
  200 & 80 & 0.0154 & 0.0657 & 0.0634 & 0.0629 & 0.0017 & 0.0416 & 0.0218 & 0.0080 & 0.0026 \\ 
  200 & 100 & 0.0059 & 0.0412 & 0.0517 & 0.0512 & 0.0016 & 0.0002 & 0.0001 & 0.0002 & 0.0001 \\ 
  500 & 20 & 0.0332 & 0.1587 & 0.1894 & 0.1908 & 0.1908 & 0.1475 & 0.1268 & 0.0436 & 0.0087 \\ 
  500 & 40 & 0.0189 & 0.1155 & 0.1343 & 0.1349 & 0.1347 & 0.0770 & 0.0621 & 0.0227 & 0.0069 \\ 
  500 & 60 & 0.0134 & 0.0889 & 0.1095 & 0.1102 & 0.1055 & 0.0412 & 0.0312 & 0.0123 & 0.0038 \\ 
  500 & 80 & 0.0084 & 0.0656 & 0.0946 & 0.0954 & 0.0826 & 0.0300 & 0.0154 & 0.0061 & 0.0021 \\ 
  500 & 100 & 0.0050 & 0.0499 & 0.0847 & 0.0853 & 0.0693 & 0.0249 & 0.0003 & 0.0001 & 0.0001 \\ 
  1000 & 20 & 0.0211 & 0.1200 & 0.1344 & 0.1349 & 0.1350 & 0.1043 & 0.0970 & 0.0471 & 0.0097 \\ 
  1000 & 40 & 0.0122 & 0.0863 & 0.0951 & 0.0954 & 0.0954 & 0.0542 & 0.0422 & 0.0199 & 0.0059 \\ 
  1000 & 60 & 0.0073 & 0.0662 & 0.0776 & 0.0779 & 0.0779 & 0.0414 & 0.0205 & 0.0098 & 0.0037 \\ 
  1000 & 80 & 0.0045 & 0.0539 & 0.0671 & 0.0675 & 0.0675 & 0.0336 & 0.0103 & 0.0047 & 0.0018 \\ 
  1000 & 100 & 0.0040 & 0.0475 & 0.0600 & 0.0603 & 0.0604 & 0.0296 & 0.0066 & 0.0005 & 0.0003 \\ 
\end{tabular}
\caption{In this table we demonstrate how feasible (orthonormal) the solutions are $G$ computed by  Algorithms \ref{alg:Ding}, \ref{alg:Mirzal} and \ref{alg:PG1} on UNION data set, i.e., in this table we report the average infeasibility of the solutions underlying Table \ref{tab:RSE_vs_beta}. 
} 
\label{tab:infeas_vs_beta}
\end{table}


Results from Table \ref{tab:infeas_vs_beta}, corresponding to $n=100,500,1000$ are depicted on Figures 
\ref{fig:rse_infeas_union_c}--\ref{fig:rse_infeas_union_j}.

\subsection{Numerical results for bi-orthonormal data (BION) }

In this subsection we provide the same type of results as in the previous subsection, but for the BION dataset.

We used almost the same setting as for UNION dataset: $\varepsilon=10^{-10}$,  maxit $=1000$, $\sigma=0.001$ and time limit =  $3600s$. Parameters $\gamma,\tau$ were slightly changed (based on experimental observations): $\gamma=0.75$ and $\tau=0.5$.
Additionally, we decided to take the same values for $\alpha,\beta$ in Algorithms \ref{alg:Mirzal}  and  \ref{alg:PG1}, since the  matrices $R$ in BION dataset are symmetric and both  orthogonality constraints are equally important.
We computed the results for values of $\alpha=\beta$ from $\{1,10,100,1000\}$.
In Tables \ref{tab:RSE_vs_beta_BINOM}--\ref{tab:infeas_vs_beta_BION}
we report average $\text{RSE}$ and average infeasibility, respectively, of the solutions obtained by Algorithms \ref{alg:Ding}, \ref{alg:Mirzal}, and \ref{alg:PG1}.
Since for this dataset we need to monitor how orthonormal are both matrices $G$ and $H$, we adapt 
the measure for infeasibility as follows: 
\bea \label{def:vio2}
\mbox{infeas}_{G,H}:= \frac{\norm{G^TG-I}_F+\norm{HH^T-I}_F}{1+\norm{I}_F}.
\eea


\begin{table}[ht]
\centering\scriptsize
\begin{tabular}{r|r||r|rrrr|rrrr}
 \multicolumn{1}{c|}{\multirow{2}{*}{$n$}} & \multicolumn{1}{c||}{$p$} & $\text{RSE}$ of  & \multicolumn{4}{c|}{$\text{RSE}$ of Alg.  \ref{alg:Mirzal}} & \multicolumn{4}{c}{$\text{RSE}$ of Alg.  \ref{alg:PG1}}\\
& \multicolumn{1}{c||}{(\% of $k$)} & Alg.  \ref{alg:Ding}  &$\beta=1$ & $\beta=10$& $\beta=100$& $\beta=1000$& $\beta=1$  & $\beta=10$& $\beta=100$& $\beta=1000$
\\   \hline
  50 & 20 & 0.7053 & 0.7053 & 0.7053 & 0.7053 & 0.8283 & 0.7053 & 0.7053 & 0.7055 & 0.8259 \\ 
  50 & 40 & 0.6108 & 0.6108 & 0.6108 & 0.6108 & 0.9066 & 0.6108 & 0.6108 & 0.6108 & 0.6631 \\ 
  50 & 60 & 0.4987 & 0.4987 & 0.4987 & 0.5442 & 0.9665 & 0.4987 & 0.4987 & 0.4987 & 0.5000 \\ 
  50 & 80 & 0.3526 & 0.3671 & 0.3742 & 0.4497 & 1.0282 & 0.3526 & 0.3796 & 0.3527 & 0.4374 \\ 
  50 & 100 & 0.0607 & 0.1712 & 0.2786 & 0.5198 & 1.0781 & 0.1145 & 0.1820 & 0.2604 & 0.3689 \\ 
  100 & 20 & 0.7516 & 0.7516 & 0.7516 & 0.7517 & 0.9070 & 0.7516 & 0.7516 & 0.7517 & 0.8224 \\ 
  100 & 40 & 0.6509 & 0.6509 & 0.6509 & 0.7174 & 0.9779 & 0.6509 & 0.6509 & 0.6509 & 0.6514 \\ 
  100 & 60 & 0.5315 & 0.5315 & 0.5315 & 0.5504 & 1.0401 & 0.5315 & 0.5315 & 0.5315 & 0.5352 \\ 
  100 & 80 & 0.3758 & 0.3787 & 0.4106 & 0.4542 & 1.1082 & 0.3801 & 0.3888 & 0.3917 & 0.3898 \\ 
  100 & 100 & 0.1377 & 0.1993 & 0.3311 & 0.4898 & 1.1734 & 0.0457 & 0.1016 & 0.2758 & 0.3757 \\ 
  200 & 20 & 0.7884 & 0.7884 & 0.7884 & 0.7884 & 0.9499 & 0.7884 & 0.7884 & 0.7884 & 0.7888 \\ 
  200 & 40 & 0.6828 & 0.6828 & 0.6828 & 0.6828 & 1.0325 & 0.6828 & 0.6828 & 0.6828 & 0.6828 \\ 
  200 & 60 & 0.5575 & 0.5575 & 0.5575 & 0.5647 & 1.0938 & 0.5575 & 0.5575 & 0.5575 & 0.5610 \\ 
  200 & 80 & 0.3942 & 0.3942 & 0.3965 & 0.5019 & 1.1618 & 0.3942 & 0.3942 & 0.3942 & 0.4373 \\ 
  200 & 100 & 0.1447 & 0.1851 & 0.3014 & 0.5400 & 1.2297 & 0.0202 & 0.1429 & 0.2964 & 0.3315 \\ 
  500 & 20 & 0.8242 & 0.8242 & 0.8242 & 0.8242 & 0.9956 & 0.8242 & 0.8242 & 0.8242 & 0.8243 \\ 
  500 & 40 & 0.7138 & 0.7138 & 0.7138 & 0.7138 & 1.0679 & 0.7138 & 0.7138 & 0.7138 & 0.7138 \\ 
  500 & 60 & 0.5828 & 0.5828 & 0.5828 & 0.6045 & 1.1534 & 0.5828 & 0.5828 & 0.5828 & 0.5828 \\ 
  500 & 80 & 0.4121 & 0.4121 & 0.4203 & 0.5285 & 1.2160 & 0.4121 & 0.4121 & 0.4121 & 0.4334 \\ 
  500 & 100 & 0.1405 & 0.1814 & 0.3401 & 0.5854 & 1.2822 & 0.0067 & 0.1059 & 0.2044 & 0.3378 \\ 
  1000 & 20 & 0.8436 & 0.8436 & 0.8436 & 0.8436 & 1.0261 & 0.8436 & 0.8436 & 0.8436 & 0.8436 \\ 
  1000 & 40 & 0.7306 & 0.7306 & 0.7306 & 0.7309 & 1.0916 & 0.7306 & 0.7306 & 0.7306 & 0.7306 \\ 
  1000 & 60 & 0.5965 & 0.5965 & 0.5965 & 0.6121 & 1.1669 & 0.5965 & 0.5965 & 0.5965 & 0.5968 \\ 
  1000 & 80 & 0.4218 & 0.4218 & 0.4256 & 0.5338 & 1.2389 & 0.4218 & 0.4218 & 0.4218 & 0.4397 \\ 
  1000 & 100 & 0.1346 & 0.1635 & 0.3324 & 0.5755 & 1.3080 & 0.0096 & 0.0697 & 0.1661 & 0.2188 \\ 
\end{tabular}
\caption{$\text{RSE}$ obtained by Algorithms \ref{alg:Ding}, \ref{alg:Mirzal} and \ref{alg:PG1} on the BION  data. For the latter two algorithms, we used  $\alpha=\beta\in \{1,10,100,1000\}$. 
For each $n\in\{50,100,200,500,1000\}$ we take all ten matrices $R$ (five of them corresponding to $k=0.2n$ and five to $k=0.4n$). We run all  three algorithms on these matrices with inner dimensions $p\in \{0.2k,0.4k,\ldots,1.0k\}$
with all possible values of $\alpha=\beta$. 
Like before, each row represents the average (arithmetic mean value) of $\text{RSE}$ obtained on instances corresponding 
to given $n$ and given $p$ as a percentage of $k$. We can see that the larger the $\beta$, the worse the $\text{RSE}$, which is consistent with expectations.
}
\label{tab:RSE_vs_beta_BINOM} 
\end{table}

\begin{table}[ht!]
\centering \scriptsize 
\begin{tabular}{r|r||r|rrrr|rrrr} 
 \multicolumn{1}{c|}{\multirow{2}{*}{$n$}} & \multicolumn{1}{c||}{$p$} & Infeas. of  & \multicolumn{4}{c|}{Infeas. of Alg.  \ref{alg:Mirzal}} & \multicolumn{4}{c}{Infeas. of Alg.  \ref{alg:PG1}}\\
& \multicolumn{1}{c||}{(\% of $k$)} & Alg.  \ref{alg:Ding}  & $\beta=1$ & $\beta=10$& $\beta=100$& $\beta=1000$ & $\beta=1$ & $\beta=10$& $\beta=100$& $\beta=1000$
\\ 
  \hline
  50 & 20 & 0.0001 & 0.0070 & 0.0036 & 0.0010 & 0.0068 & 0.0017 & 0.0021 & 0.0021 & 0.0026 \\ 
  50 & 40 & 0.0000 & 0.0041 & 0.0021 & 0.0004 & 0.0056 & 0.0008 & 0.0012 & 0.0012 & 0.0014 \\ 
  50 & 60 & 0.0000 & 0.0030 & 0.0009 & 0.0032 & 0.0038 & 0.0005 & 0.0008 & 0.0009 & 0.0009 \\ 
  50 & 80 & 0.0000 & 0.0183 & 0.0030 & 0.0021 & 0.0028 & 0.0004 & 0.0202 & 0.0006 & 0.0013 \\ 
  50 & 100 & 0.0355 & 0.0533 & 0.0127 & 0.0045 & 0.0027 & 0.0418 & 0.0478 & 0.0123 & 0.0021 \\ 
  100 & 20 & 0.0001 & 0.0051 & 0.0024 & 0.0006 & 0.0063 & 0.0010 & 0.0012 & 0.0013 & 0.0016 \\ 
  100 & 40 & 0.0000 & 0.0029 & 0.0017 & 0.0066 & 0.0040 & 0.0004 & 0.0006 & 0.0007 & 0.0007 \\ 
  100 & 60 & 0.0000 & 0.0019 & 0.0008 & 0.0009 & 0.0027 & 0.0003 & 0.0004 & 0.0005 & 0.0005 \\ 
  100 & 80 & 0.0000 & 0.0039 & 0.0048 & 0.0015 & 0.0021 & 0.0062 & 0.0149 & 0.0037 & 0.0006 \\ 
  100 & 100 & 0.0606 & 0.0454 & 0.0105 & 0.0022 & 0.0018 & 0.0106 & 0.0228 & 0.0173 & 0.0028 \\ 
  200 & 20 & 0.0002 & 0.0033 & 0.0019 & 0.0005 & 0.0043 & 0.0005 & 0.0007 & 0.0007 & 0.0007 \\ 
  200 & 40 & 0.0001 & 0.0017 & 0.0010 & 0.0002 & 0.0027 & 0.0002 & 0.0003 & 0.0004 & 0.0003 \\ 
  200 & 60 & 0.0001 & 0.0010 & 0.0005 & 0.0004 & 0.0019 & 0.0001 & 0.0002 & 0.0002 & 0.0004 \\ 
  200 & 80 & 0.0000 & 0.0006 & 0.0006 & 0.0015 & 0.0014 & 0.0001 & 0.0001 & 0.0002 & 0.0013 \\ 
  200 & 100 & 0.0425 & 0.0280 & 0.0064 & 0.0019 & 0.0015 & 0.0046 & 0.0224 & 0.0240 & 0.0034 \\ 
  500 & 20 & 0.0001 & 0.0017 & 0.0011 & 0.0003 & 0.0025 & 0.0002 & 0.0003 & 0.0003 & 0.0003 \\ 
  500 & 40 & 0.0001 & 0.0008 & 0.0005 & 0.0001 & 0.0016 & 0.0001 & 0.0001 & 0.0002 & 0.0002 \\ 
  500 & 60 & 0.0000 & 0.0005 & 0.0003 & 0.0006 & 0.0013 & 0.0001 & 0.0001 & 0.0001 & 0.0002 \\ 
  500 & 80 & 0.0000 & 0.0003 & 0.0009 & 0.0009 & 0.0008 & 0.0000 & 0.0001 & 0.0001 & 0.0016 \\ 
  500 & 100 & 0.0258 & 0.0184 & 0.0045 & 0.0013 & 0.0007 & 0.0017 & 0.0101 & 0.0175 & 0.0053 \\ 
  1000 & 20 & 0.0001 & 0.0010 & 0.0006 & 0.0002 & 0.0024 & 0.0001 & 0.0002 & 0.0002 & 0.0002 \\ 
  1000 & 40 & 0.0000 & 0.0005 & 0.0003 & 0.0001 & 0.0009 & 0.0001 & 0.0002 & 0.0003 & 0.0002 \\ 
  1000 & 60 & 0.0000 & 0.0003 & 0.0002 & 0.0004 & 0.0009 & 0.0003 & 0.0002 & 0.0003 & 0.0003 \\ 
  1000 & 80 & 0.0000 & 0.0002 & 0.0005 & 0.0007 & 0.0006 & 0.0040 & 0.0001 & 0.0002 & 0.0020 \\ 
  1000 & 100 & 0.0173 & 0.0117 & 0.0031 & 0.0009 & 0.0005 & 0.0043 & 0.0050 & 0.0121 & 0.0060 
\end{tabular}
\caption{In this table we demonstrate how feasible (orthonormal) are the solutions $G$ and $H$ computed by Algorithms \ref{alg:Ding}, \ref{alg:Mirzal}, and \ref{alg:PG1} on the BION dataset, i.e., in this table we report the average infeasibility \eqref{def:vio2} of the solutions underlying Table \ref{tab:RSE_vs_beta_BINOM}. 
We can observe that with these settings of all algorithms we can bring infeasibility to order of $10^{-3}$ very often, for all values of $\beta$.
} 
\label{tab:infeas_vs_beta_BION}
\end{table}

\begin{figure}
\centering
  \begin{subfigure}{6cm} 
    \centering\includegraphics[width=5cm]{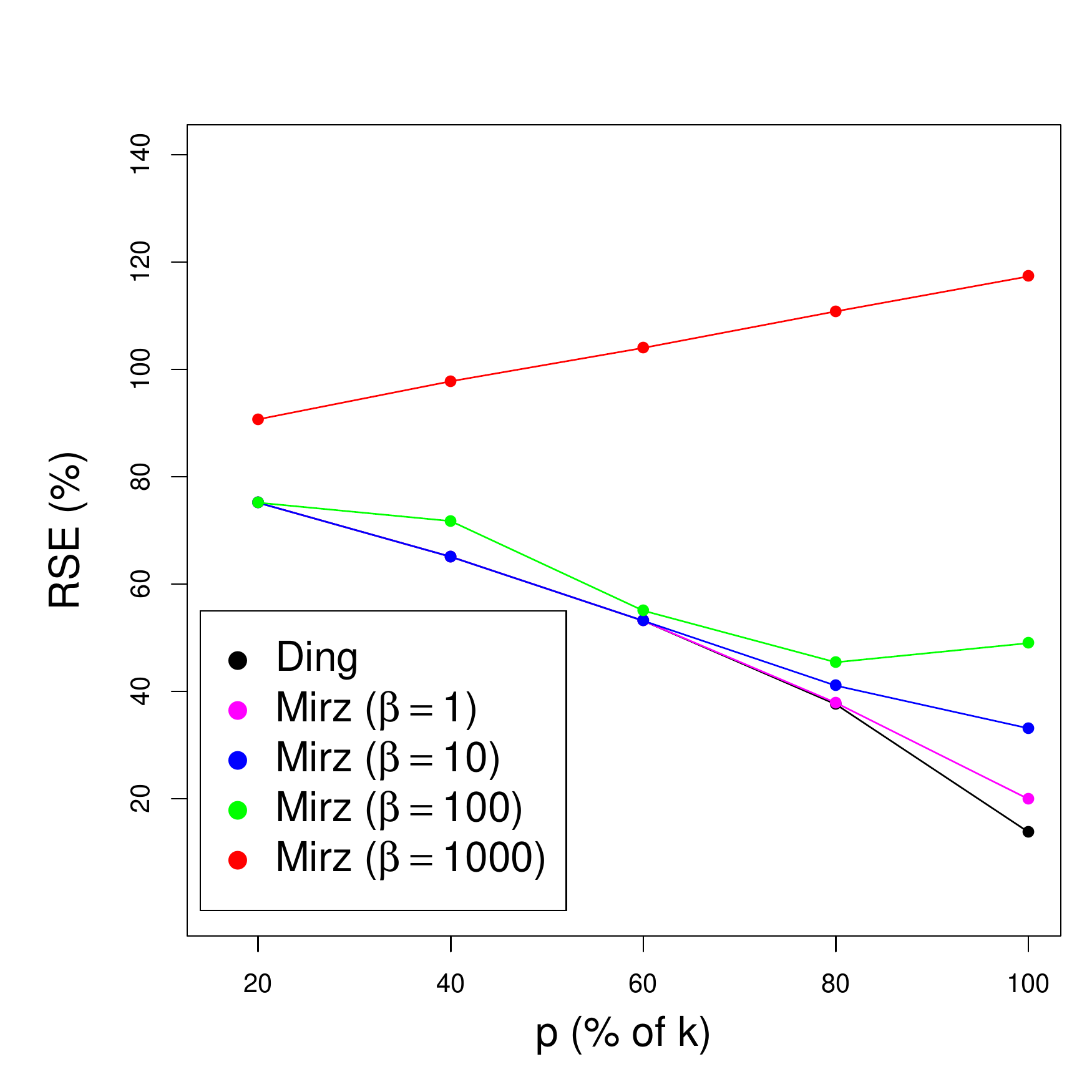}
    \caption{Values of $\text{RSE}$ for different values of $\beta$ obtained by  Algorithms \ref{alg:Ding} and \ref{alg:Mirzal} in BION data with $n=100$}\label{fig:rse_bion_a}
  \end{subfigure}
   \begin{subfigure}{0.7cm}
  ~
  \end{subfigure}
  \begin{subfigure}{6cm}
    \centering\includegraphics[width=5cm]{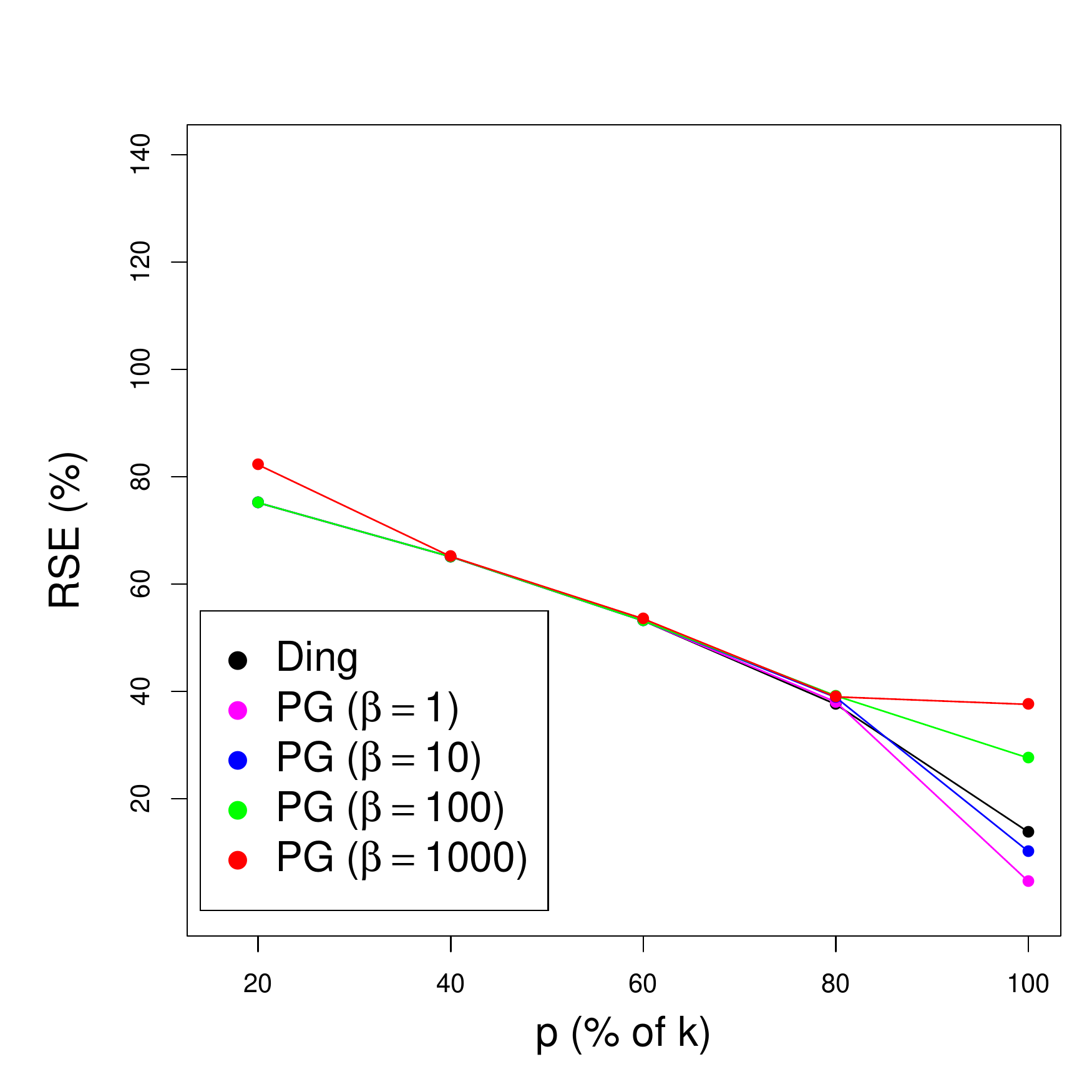}
    \caption{Values of $\text{RSE}$ for different values of $\beta$ obtained by  Algorithms \ref{alg:Ding} and  \ref{alg:PG1} on BION data with $n=100$.}\label{fig:rse_bion_b}  
    \end{subfigure}\\
    
  \begin{subfigure}{6cm}
    \centering\includegraphics[width=5cm]{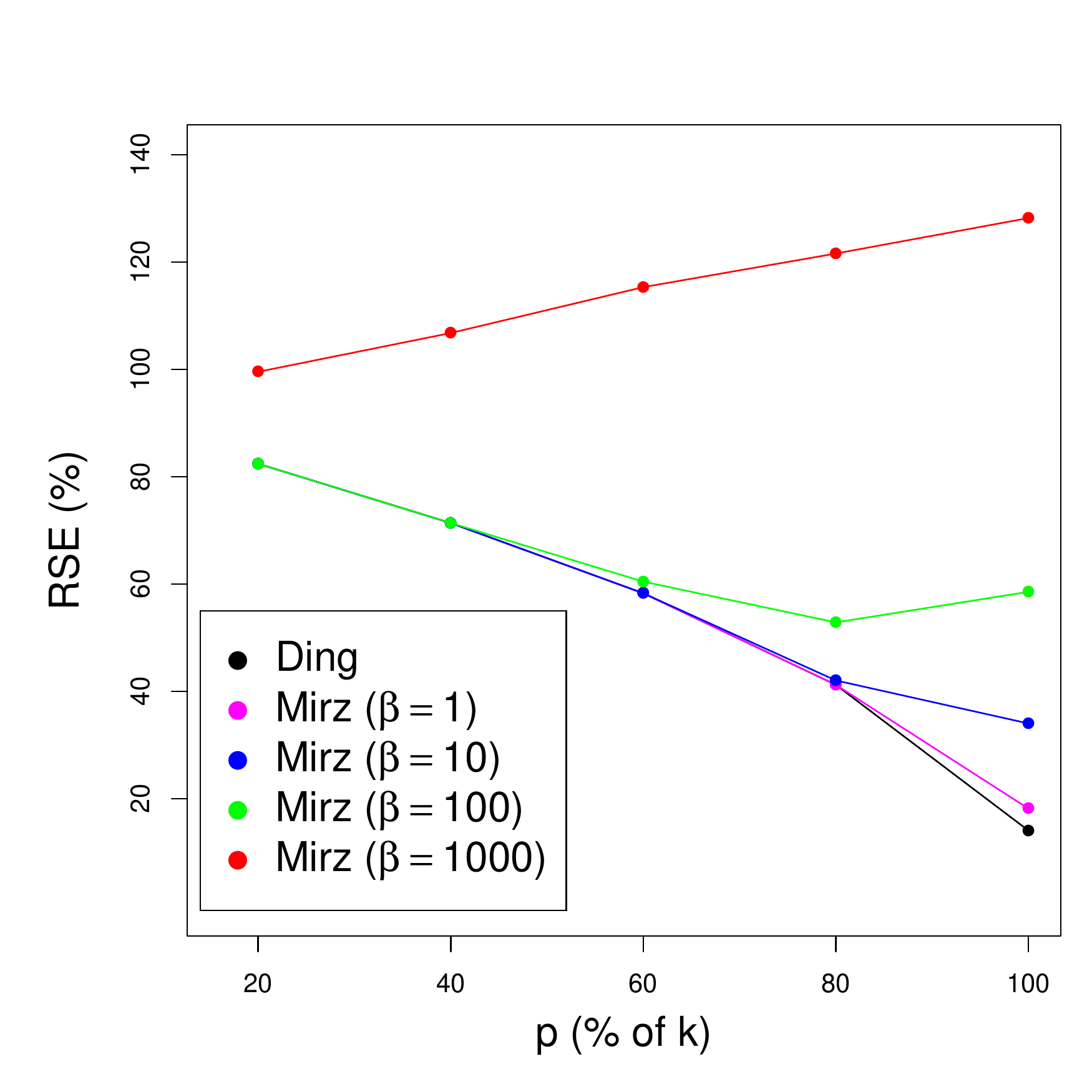}
    \caption{Values of $\text{RSE}$ for different values of $\beta$ obtained by  Algorithms \ref{alg:Ding} and \ref{alg:Mirzal}  in BION data with $n=500$}\label{fig:rse_bion_c}
  \end{subfigure}
   \begin{subfigure}{0.7cm}
  ~
  \end{subfigure}
  \begin{subfigure}{6cm}
    \centering\includegraphics[width=5cm]{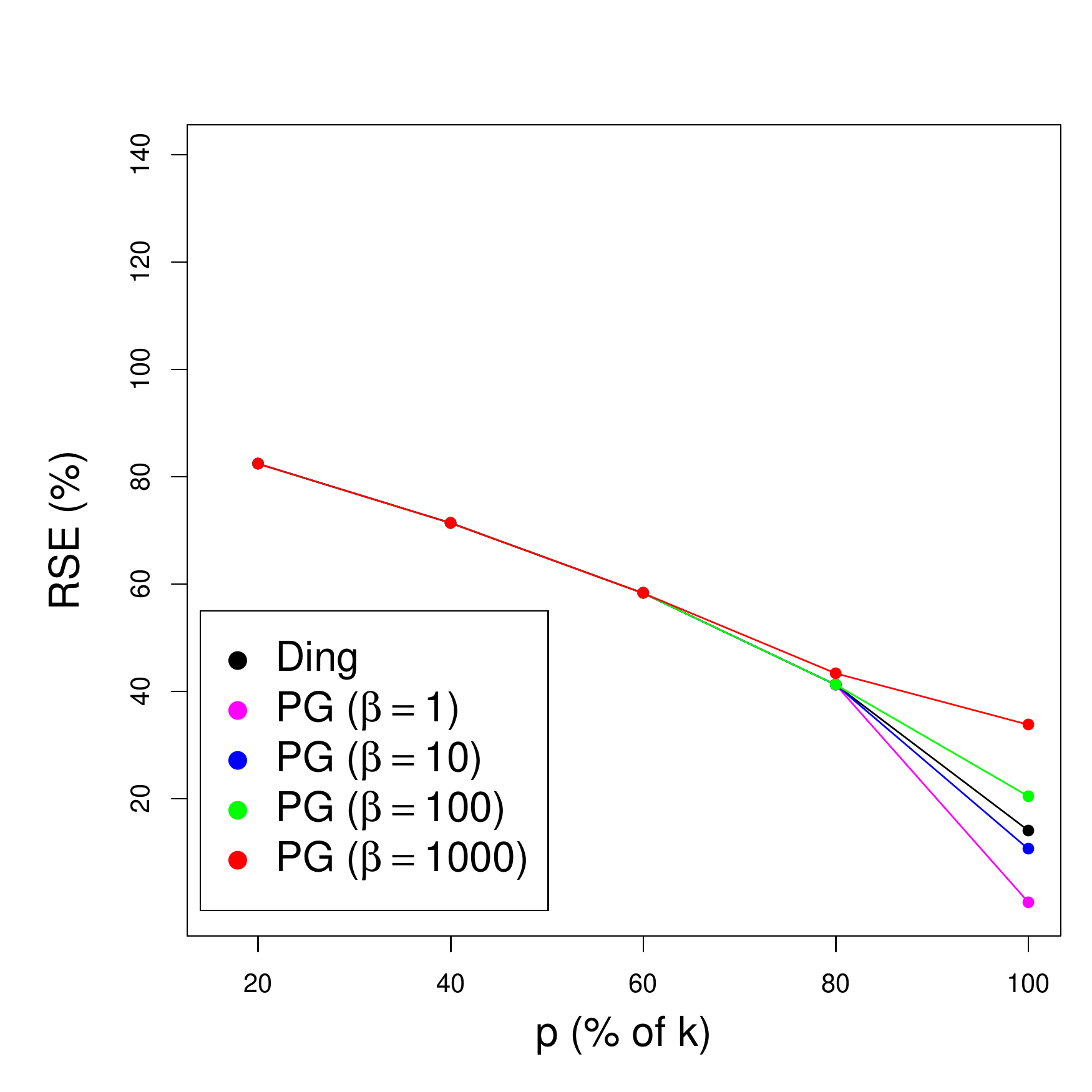}
    \caption{Values of $\text{RSE}$ for different values of $\beta$ obtained by  Algorithms \ref{alg:Ding} and  \ref{alg:PG1} on BION data with $n=500$.}\label{fig:rse_bion_d}
  \end{subfigure}\\
  \begin{subfigure}{6cm}
    \centering\includegraphics[width=5cm]{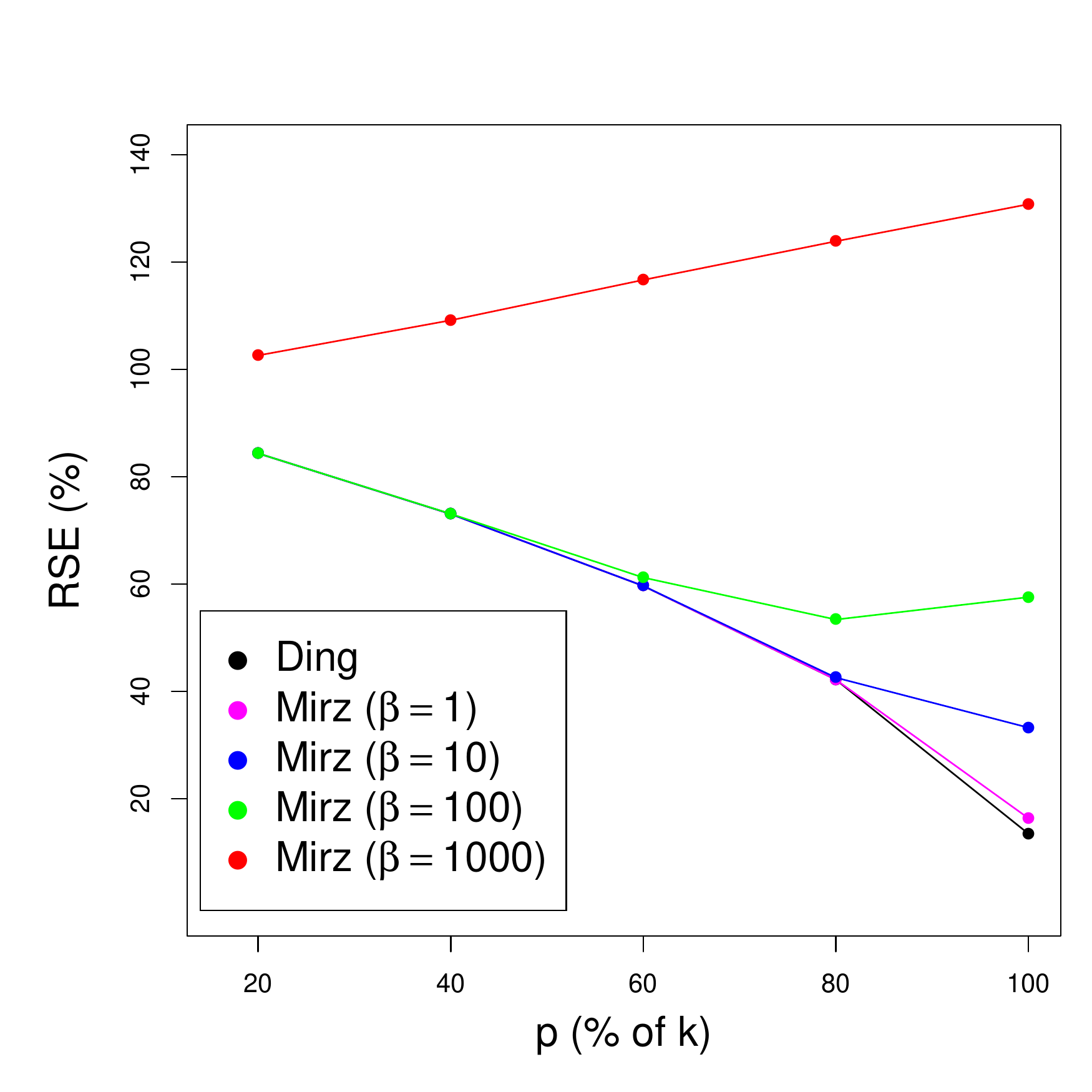}
    \caption{Values of $\text{RSE}$ for different values of $\beta$ obtained by  Algorithms \ref{alg:Ding} and  \ref{alg:Mirzal} in BION data with $n=1000$}\label{fig:rse_bion_e}
  \end{subfigure}
   \begin{subfigure}{0.7cm}
  ~
  \end{subfigure}
  \begin{subfigure}{6cm}
    \centering\includegraphics[width=5cm]{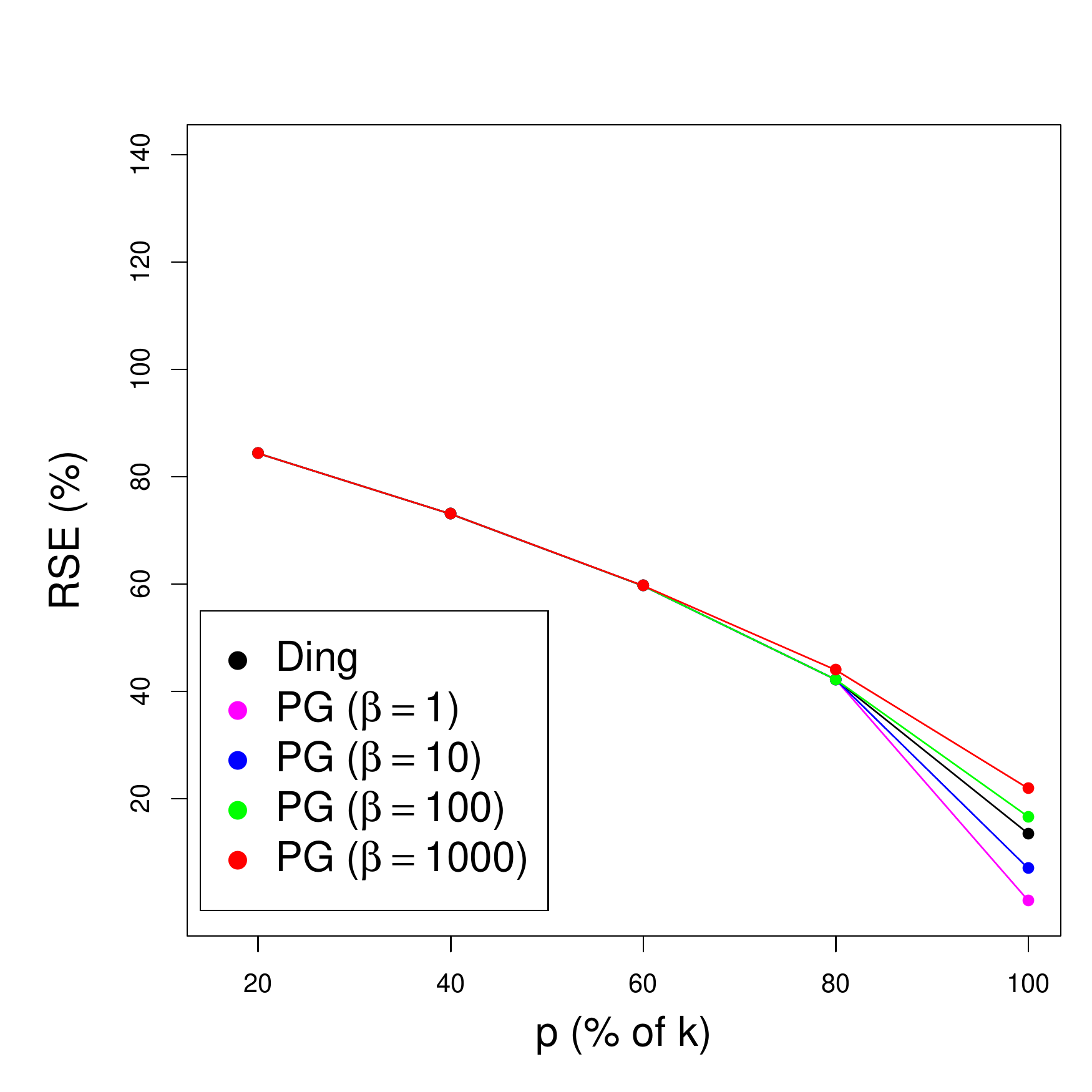}
    \caption{Values of $\text{RSE}$ for different values of $\beta$ obtained by  Algorithms \ref{alg:Ding} and  \ref{alg:PG1} on BION data with $n=1000$.}\label{fig:rse_bion_f}
    \end{subfigure}
    \caption{This figure contains six plots which illustrate the quality of Algorithms \ref{alg:Ding}, \ref{alg:Mirzal} and \ref{alg:PG1} regarding $\text{RSE}$ on BION instances with  $n=100,500,1000$ and  $k=0.2n,~0.4n$, for  $\beta\in\{1,10,100,1000\}$. We can observe that Algorithm \ref{alg:PG1} is more stable, less dependent to the choice of $\beta$ and 
    is computing better values of  RSE.  }\label{fig:rse_bion}
\end{figure}
\begin{figure}
\centering
  \begin{subfigure}{6cm} 
    \centering\includegraphics[width=5cm]{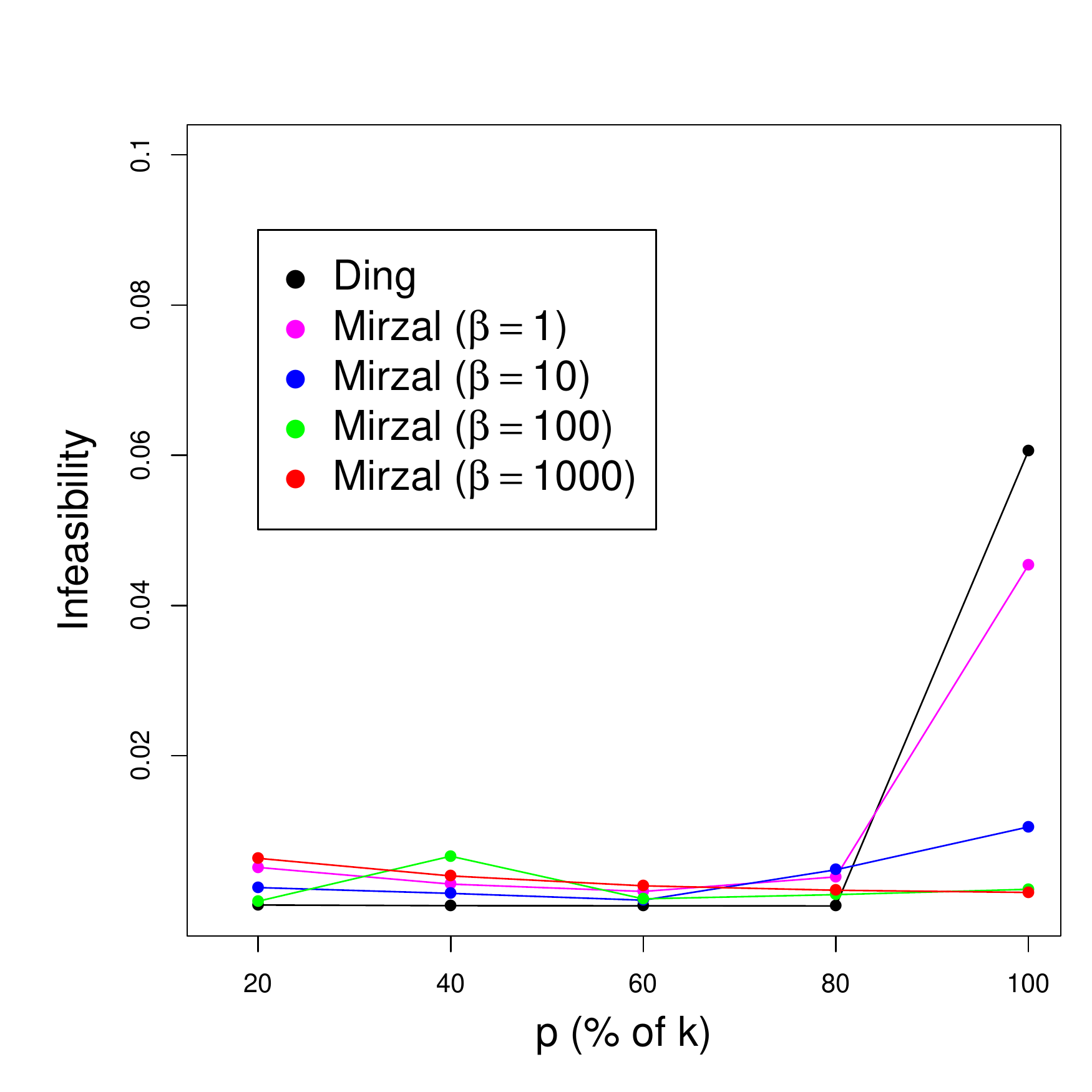}
    \caption{Values of infeasibility  for different values of $\beta$ obtained by  Algorithms \ref{alg:Ding} and  \ref{alg:Mirzal} on BION data with $n=100$}\label{fig:rse_infeas_bion_a}
  \end{subfigure}
   \begin{subfigure}{0.7cm}
  ~
  \end{subfigure}
  \begin{subfigure}{6cm}
    \centering\includegraphics[width=5cm]{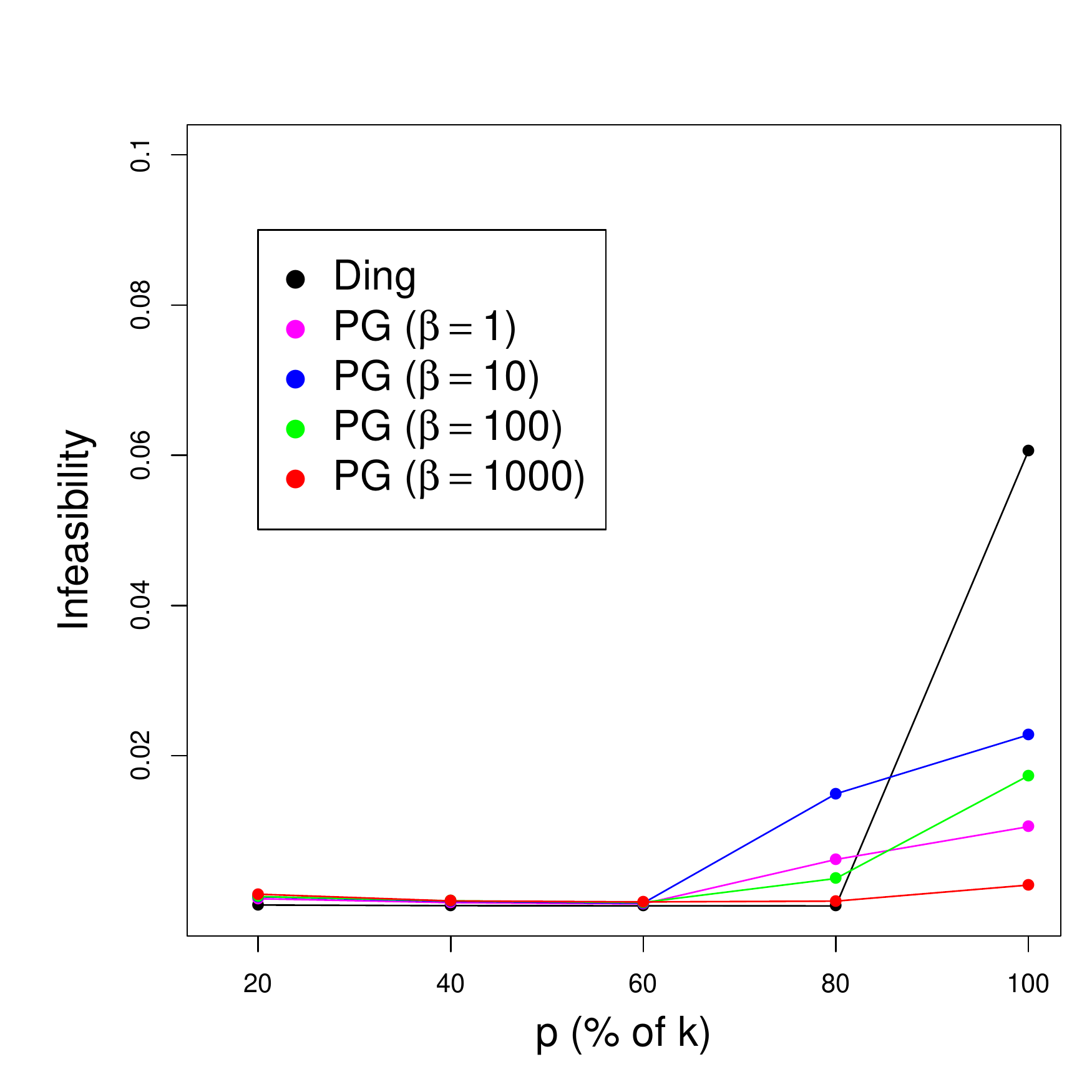}
    \caption{Values of infeasibility for different values of $\beta$ obtained by  Algorithms \ref{alg:Ding} and  \ref{alg:PG1} on BION data with $n=100$.}\label{fig:rse_infeas_bion_b}  
    \end{subfigure}
    
  \begin{subfigure}{6cm}
    \centering\includegraphics[width=5cm]{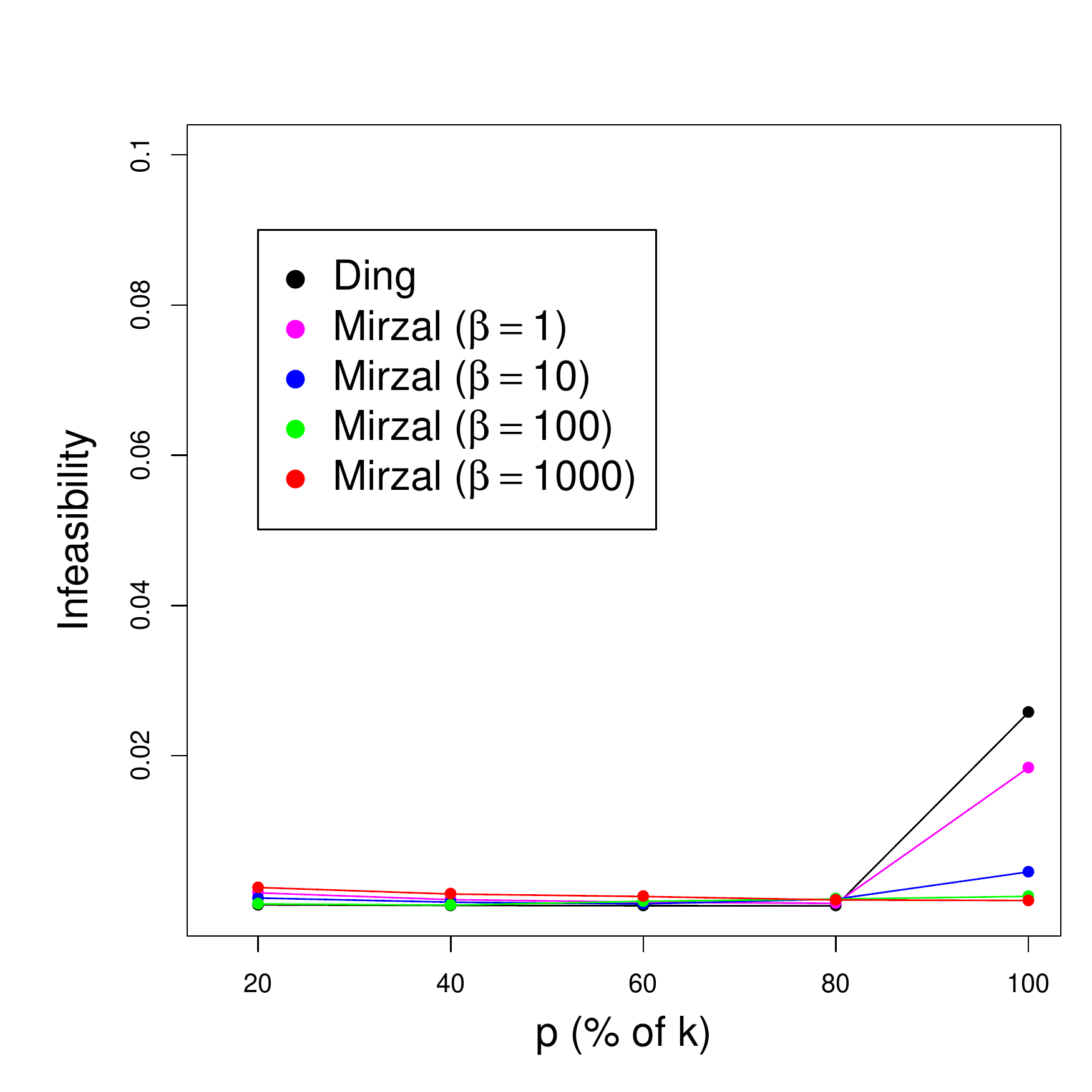}
    \caption{Values of infeasibility for different values of $\beta$ obtained by  Algorithms \ref{alg:Ding} and  \ref{alg:Mirzal} in BION data with $n=500$}\label{fig:rse_infeas_bion_c}
  \end{subfigure}
   \begin{subfigure}{0.7cm}
  ~
  \end{subfigure}
  \begin{subfigure}{6cm}
    \centering\includegraphics[width=5cm]{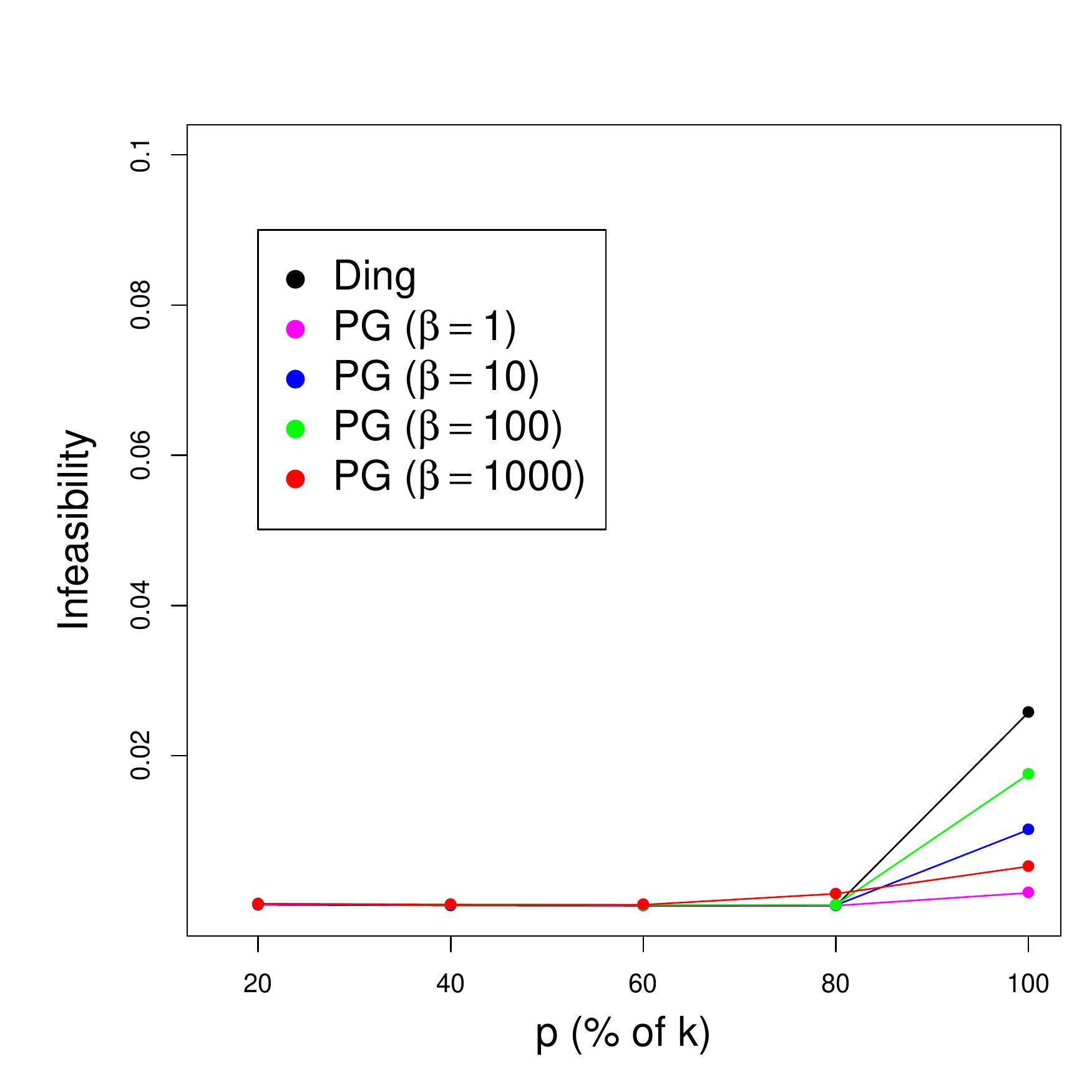}
    \caption{Values of infeasibility for different values of $\beta$ obtained by  Algorithms \ref{alg:Ding} and  \ref{alg:PG1} on BION data with $n=500$.}\label{fig:rse_infeas_bion_d}
    \end{subfigure}
    
      \begin{subfigure}{6cm}
    \centering\includegraphics[width=5cm]{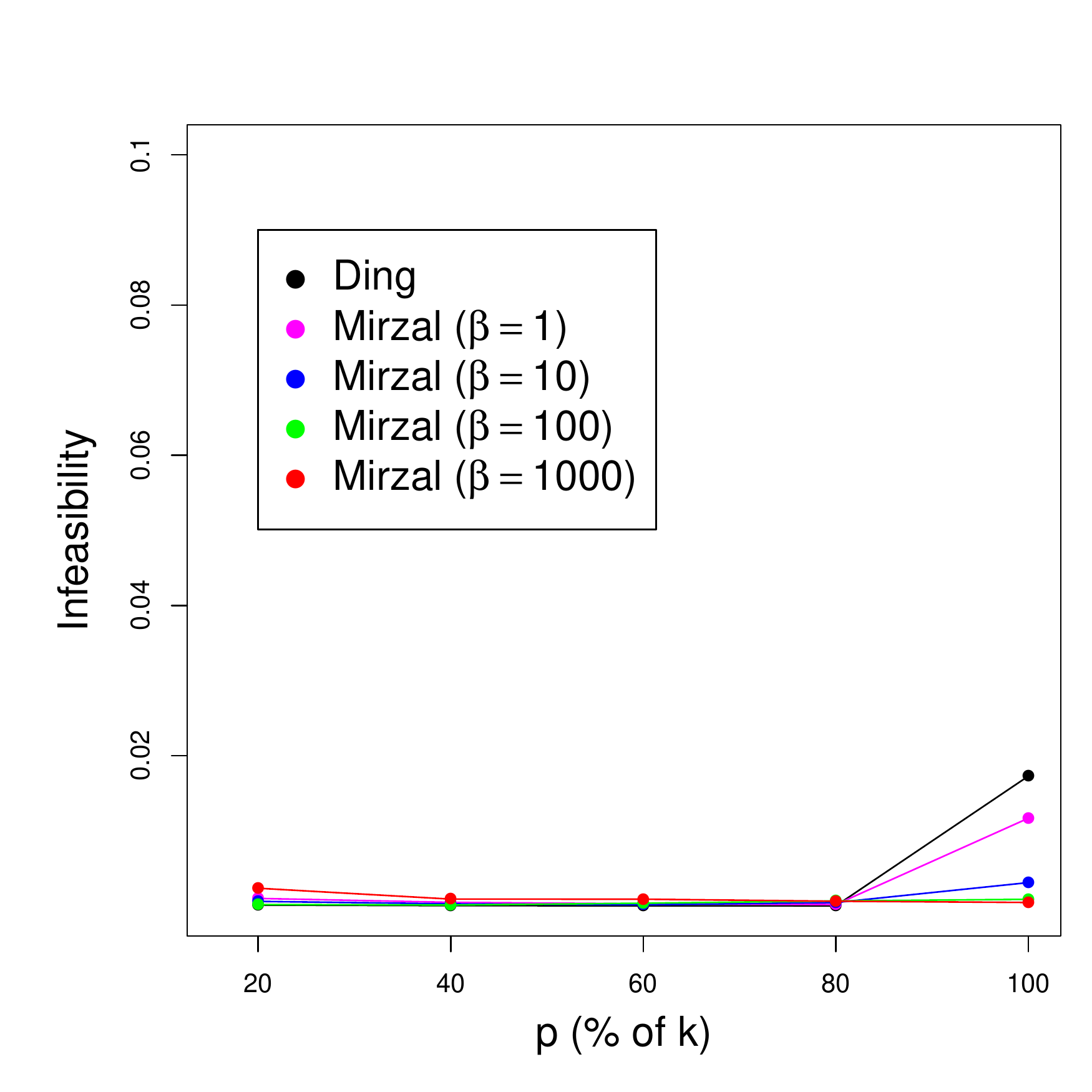}
    \caption{Values of infeasibility for different values of $\beta$ obtained by  Algorithms \ref{alg:Ding} and  \ref{alg:Mirzal} in BION data with $n=1000$}\label{fig:rse_infeas_bion_e}
  \end{subfigure}
   \begin{subfigure}{0.7cm}
  ~
  \end{subfigure}
  \begin{subfigure}{6cm}
    \centering\includegraphics[width=5cm]{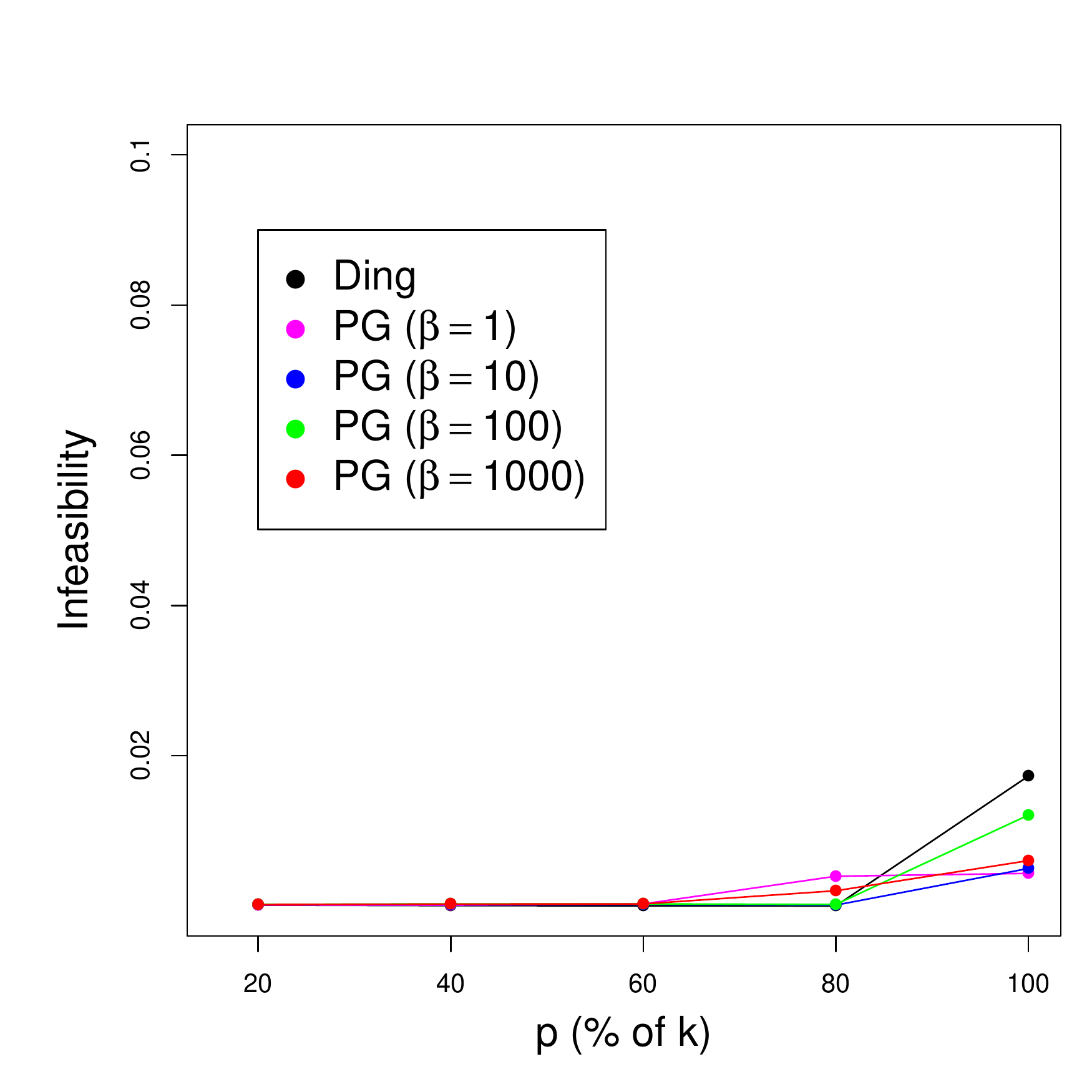}
    \caption{Values of infeasibility for different values of $\beta$ obtained by  Algorithms \ref{alg:Ding} and \ref{alg:PG1} on BION data with $n=1000$.}\label{fig:rse_infeas_bion_f}
    \end{subfigure}
    \caption{This figure contains six plots which illustrate the quality of Algorithms \ref{alg:Ding}, \ref{alg:Mirzal} and \ref{alg:PG1} regarding the infeasibility on BION instances with  $n=100,500,1000$ and  $k=0.2n,~0.4n$,  for  $\beta\in\{1,10,100,1000\}$. We can observe that Algorithm \ref{alg:PG1} computes solutions with infeasibility \eqref{def:vio2} slightly smaller compared to solutions computed by Algorithm \ref{alg:Mirzal}. }\label{fig:rse_infeas_bion}
\end{figure}

Figures \ref{fig:rse_bion_a}--\ref{fig:rse_bion_f} and \ref{fig:rse_infeas_bion_a}--\ref{fig:rse_infeas_bion_f} depict $\text{RSE}$ and infeasibility reached by the three compared algorithms, for $n=100,500,1000$. We can see that all three algorithms behave well, however, Algorithm \ref{alg:PG1} is more stable and less dependent on the choice of $\beta$. 
It is interesting to see that $\beta$ does not have a big impact on RSE and infeasibility for  Algorithm \ref{alg:PG1}, a significant difference can be observed only when the internal dimension is equal to the real internal dimension, i.e., when $p=100 \%$. Based on these numerical results, we can conclude that smaller $\beta$ achieve better RSE and almost the same infeasibility, so it would make sense to use $\beta=1$.

For Algorithm \ref{alg:Mirzal} these differences are bigger and it is less obvious which $\beta$ is appropriate. Again, if RSE is more important then smaller values of $\beta$ should be taken, otherwise larger values. 


\section{Concluding remarks}

We presented a projected gradient method to solve the orthogonal non-negative matrix factorization problem. We penalized the deviation from orthonormality with some positive parameters and added the resulted terms to the objective function of the standard non-negative matrix factorization problem. Then, we considered minimizing the resulted objective function under the non-negativity conditions only, in a block coordinate decent approach. 

The method was tested on two sets of synthetic data, one containing uni-orthonormal matrices and the other containing bi-orthonormal matrices. Different values for the adjusting parameters of orthogonality were applied in the implementation to determine good pairs of such values. The performance of our algorithm was compared with two algorithms based on  multiplicative updates rules. Algorithms were compared regarding the quality of factorization (RSE) and how much the resulting factors deviate from orthonormality.

We provided an extensive list of numerical results which demonstrate that our method is very competitive and outperforms the others. 
 
 \section*{Acknowledgment}
The work of the first author is supported by the Swiss Government Excellence Scholarships grant number ESKAS-2019.0147. This author also thanks the University of Applied Sciences and Arts, Northwestern Switzerland for supporting the work.

The work of the second author was partially funded by Slovenian Research Agency under research program P2-0256 and research projects N1-0057, N1-0071, and J1-8155.

The authors would also like to thank to Andri Mirzal (Faculty of Computing, Universiti Teknologi Malaysia) for providing the code for his algorithm (Algorithm \ref{alg:Mirzal}) to solve \eqref{ONMF}. This code was also adapted by the authors to solve \eqref{biONMF}.

\bigskip

\bibliographystyle{acm}
\bibliography{onmf}

\end{document}